\documentclass[a4paper,12pt]{amsart}
\usepackage{amssymb}
\usepackage{ifthen}
 \usepackage[dvips]{graphicx}
\nonstopmode \numberwithin{equation}{section}
\usepackage{amssymb}
\usepackage{ifthen}
\usepackage{graphicx}
\usepackage{amsmath}
\usepackage[T1]{fontenc} 
\usepackage[utf8]{inputenc}
\usepackage[usenames,dvipsnames]{color}
\usepackage{color}
\usepackage[english]{babel}
\usepackage{fancyhdr}
\usepackage{fancybox}
\usepackage{tikz}

\nonstopmode \numberwithin{equation}{section}
\theoremstyle{plain}
\newtheorem{thm}[equation]{Theorem}
\newtheorem{cor}[equation]{Corollary}
\newtheorem{lem}[equation]{Lemma}
\newtheorem{prop}{Proposition}

\newtheorem{conj}{Conjecture}

\theoremstyle{definition}
\newtheorem{defn}{Definition}[section]

\newtheorem{prob}{Problem}
\newtheorem{rem}{Remark}[section]

\newcounter{minutes}
\divide\time by 60
\newcounter{hours}
\multiply\time by 60
\addtocounter{minutes}{-\time}

\newcounter {own}
\def\theown {\thesection       .\arabic{own}}

\newenvironment{pf}[1][]{%
 \vskip 3mm
 \noindent
 \ifthenelse{\equal{#1}{}}%
  {{\slshape Proof. }}%
  {{\slshape #1.} }%
 }%
{\qed\bigskip}

\newcounter{alphabet}

\def\be{\begin{equation}}
\def\ee{\end{equation}}

\newcommand{\bee}{\begin{enumerate}}
\newcommand{\eee}{\end{enumerate}}

\newcommand{\blem}{\begin{lem}}
\newcommand{\elem}{\end{lem}}
\newcommand{\bthm}{\begin{thm}}
\newcommand{\ethm}{\end{thm}}
\newcommand{\bcor}{\begin{cor}}
\newcommand{\ecor}{\end{cor}}
\newcommand{\beg}{\begin{examp}}
\newcommand{\eeg}{\end{examp}}
\newcommand{\begs}{\begin{examples}}
\newcommand{\eegs}{\end{examples}}

\newcommand{\bdefn}{\begin{defn}}
\newcommand{\edefn}{\end{defn}}

\newcommand{\bprob}{\begin{prob}}
\newcommand{\eprob}{\end{prob}}
\newcommand{\bei}{\begin{itemize}}
\newcommand{\eei}{\end{itemize}}

\newcommand{\bcon}{\begin{conj}}
\newcommand{\econ}{\end{conj}}
\newcommand{\bcons}{\begin{conjs}}
\newcommand{\econs}{\end{conjs}}
\newcommand{\bprop}{\begin{prop}}
\newcommand{\eprop}{\end{prop}}
\newcommand{\br}{\begin{rem}}
\newcommand{\er}{\end{rem}}
\newcommand{\brs}{\begin{rems}}
\newcommand{\ers}{\end{rems}}
\newcommand{\bo}{\begin{obser}}
\newcommand{\eo}{\end{obser}}
\newcommand{\bos}{\begin{obsers}}
\newcommand{\eos}{\end{obsers}}
\newcommand{\bpf}{\begin{pf}}
\newcommand{\epf}{\end{pf}}
\newcommand{\ba}{\begin{array}}
\newcommand{\ea}{\end{array}}
\newcommand{\beq}{\begin{eqnarray}}
\newcommand{\beqq}{\begin{eqnarray*}}
\newcommand{\eeq}{\end{eqnarray}}
\newcommand{\eeqq}{\end{eqnarray*}}

\begin{document}

\title{Bohr phenomenon for certain classes of harmonic mappings}

\author{Molla Basir Ahamed}
\address{Molla Basir Ahamed,
	School of Basic Sciences,
	Indian Institute of Technology Bhubaneswar,
	Bhubaneswar-752050, Odisha, India.}
\email{mba15@iitbbs.ac.in}

\author{Vasudevarao Allu}
\address{Vasudevarao Allu,
School of Basic Sciences,
Indian Institute of Technology Bhubaneswar,
Bhubaneswar-752050, Odisha, India.}
\email{avrao@iitbbs.ac.in}

\subjclass[{AMS} Subject Classification:]{Primary 30C45, 30C50, 30C80}
\keywords{Analytic, univalent, harmonic functions; starlike, convex, close-to-convex functions; coefficient estimate, growth theorem, Bohr radius.}

\def\thefootnote{}
\footnotetext{ {\tiny File:~\jobname.tex,
printed: \number\year-\number\month-\number\day,
          \thehours.\ifnum\theminutes<10{0}\fi\theminutes }
} \makeatletter\def\thefootnote{\@arabic\c@footnote}\makeatother

\begin{abstract}
Bohr phenomenon for analytic functions $ f $  where $ f(z)=\sum_{n=0}^{\infty}a_nz^n $, first introduced by Harald Bohr in $ 1914 $, deals with finding the largest radius $ r_f $, $ 0<r_f<1 $, such that the inequality $ \sum_{n=0}^{\infty}|a_nz^n|<1 $ holds whenever $ |f(z)|<1 $ holds in the unit disk $ \mathbb{D}=\{z\in\mathbb{C} : |z|<1\} $. The Bohr phenomenon for the harmonic functions of the form $ f(z)=h+\overline{g} $, where $ h(z)=\sum_{n=0}^{\infty}a_nz^n $ and $ g(z)=\sum_{n=1}^{\infty}b_nz^n $ is to find the largest radius $ r_f $, $ 0<r_f<1 $ such that
\begin{equation*}
	\sum_{n=1}^{\infty}\left(|a_n|+|b_n|\right)|z|^n\leq d(f(0),\partial f(\mathbb{D}))
\end{equation*}
holds for $ |z|\leq r_f $, where $ d(f(0),\partial f(\mathbb{D})) $ is the Euclidean distance between $ f(0) $ and the boundary of $ f(\mathbb{D}) $. In this paper, we prove several improved versions of the sharp Bohr radius for the classes of harmonic and univalent functions.  Further, we prove several corollaries as a consequence of the main results. 
\end{abstract}

\maketitle
\pagestyle{myheadings}
\markboth{Molla Basir Ahamed and Vasudevarao Allu}{Bohr phenomenon for certain classes of harmonic mappings}

\section{Introduction}

\par Let $ \mathcal{B} $ be the class of analytic functions $ f $ in the unit disk $ \mathbb{D}=\{z\in\mathbb{C} : |z|<1\} $ such that $ |f(z)|<1 $ for all $ z\in\mathbb{D} $, and let $ \mathcal{B}_0=\{f\in\mathcal{B} : f(0)=0\} $. In $ 1914 $, Bohr \cite{Bohr-1914} proved that if $ f\in\mathcal{B} $ is of the form $ f(z)=\sum_{n=0}^{\infty}a_nz^n $, then the majorant series $ M_f(r)= \sum_{n=0}^{\infty}|a_n||z|^n $ of $ f $ satisfies 
\begin{equation}\label{e-11.1}
	M_{f_0}(r)=\sum_{n=1}^{\infty}|a_n||z|^n\leq 1-|a_0|=d(f(0),\partial f(\mathbb{D}))
\end{equation} 
for all $ z\in\mathbb{D} $ with $ |z|=r\leq 1/3 $, where $ f_0(z)=f(z)-f(0) $. The interesting inequality \eqref{e-11.1} is known as the Bohr inequality. Bohr actually obtained the inequality \eqref{e-11.1} for $ |z|\leq 1/6 $. Later, Wiener, Riesz and Schur, independently established the inequaliy \eqref{e-11.1} for $ |z|\leq 1/3 $ and hence proved that $ 1/3 $ is the best possible. The radius $ 1/3 $ is called as Bohr radius for the class $ \mathcal{B} $. Moreover, for 
\begin{equation*}
	\phi_a(z)=\frac{a-z}{1-az},\;\; a\in [0,1),
\end{equation*}
it follows easily that $ M_{\phi_a}(r)>1 $ if, and only if, $ r>1/(1+2a) $, which for $ a\rightarrow 1 $ shows that $ 1/3 $ is optimal. \vspace{1.5mm}

\par The Bohr inequality has also emerged as an active area of research for operator algebraists after Dixon \cite{Dixon & BLMS & 1995} used it to settle in the negative a conjecture that a Banach algebra satisfying a non-unital von Neumann inequality is necessarily an operator algebra. Subsequently, Paulson and Singh \cite{Paulsen-PAMS-2004}, and Blasco \cite{Blasco-2010} have extended the Bohr inequality in the context of Banach algebra.\vspace{1.5mm}

\par The main aim of this paper is to establish various Bohr inequalities for the harmonic functions in a suitable fashion. A harmonic mapping in the unit disk $ \mathbb{D} $ is a complex-valued function $ f=u+iv $ of $ z=x+iy $ in $ \mathbb{D} $, which satisfies the Laplace equation $ \Delta f=4f_{z\bar{z}}=0 $, where $ f_{z}=(f_{x}-if_{y})/2 $ and $ f_{\bar{z}}=(f_{x}+if_{y})/2 $ and $ u $ and $ v $ are real-valued harmonic functions in $ \mathbb{D} $. It follows that the function $ f $ admits the canonical representation $ f=h+\bar{g} $, where $ h $ and $ g $ are analytic in $ \mathbb{D} $. The Jacobian $ J_f $ of $ f=h+\overline{g} $ is given by $ J_f=|h^{\prime}|^2-|g^{\prime}|^2 $. We say that $ f $ is sense-preserving in $ \mathbb{D} $ if $ J_f(z)>0 $ in $ \mathbb{D} $. Consequently, $ f $ is locally univalent and sense-preserving in $ \mathbb{D} $ if, and only if, $ J_f(z)>0 $ in $ \mathbb{D} $; or equivalently if $ h^{\prime}\neq 0 $ in $ \mathbb{D} $ and the dilation $ \omega_f := \omega=g^{\prime}/h^{\prime} $ has the property that $ |\omega(z)|<1 $ in $ \mathbb{D} $. \vspace{1.5mm}

\par  The interests in the Bohr phenomenon was revived in the $ 1990 $s due to the extensions to holomorphic function in $ \mathbb{C}^n $ and to abstract setting \cite{Aizenberg-PAMS-2000}. For example, Boas and Khavinshon \cite{Boas-1997} found bounds for Bohr's radius in any complete Reinhard domains and proved that the Bohr radius decreses to zero as the dimension of the domain incresses. In the recent years, studying Bohr inequalities and Bohr radius is become an active area of research in univalent function theory. The Bohr's phenomenon for the complex-valued harmonic mappings have been studied extensively by many authors (see \cite{Abu-CVEE-2010,Abu-JMAA-2014,Himadri-Vasu-P1,Nirupam-VVEE-2018,Nirupam-MonatsMath-2019}). In $ 2016 $, Ali \emph{et al.} \cite{Ali & Abdul & Ng & CVEE & 2016} studied Bohr radius for the stralike log-harmonic mappings and obtain several interesting results. In $ 2021 $, Liu and Ponnusamy \cite{Liu-Ponnusamy-PAMS-2021} established a version of multidimensional analogue of the refined Bohr inequality as well as improved Bohr inequality with initial coefficient being zero. Bohr-type inequalities for the class of harmonic $ p$-symmetric mappings and also for harmonic mappings with a multiple zero at the origin have been discussed by Huang \emph{et al.} \cite{Huang-Liu-Ponnusamy-MJM-2021}.  The Bohr radius for various classes of functions, for example, locally univalent harmonic mappings, $ k $-quasiconformal mappings, bounded harmonic functions, lacunary series have been studied extensively in \cite{Ismagilov-2020-JMAA,Kayumov-2018-JMAA,Kay & Pon & Sha & MN & 2018}. For more intriguing aspects of Bohr phenomenon we refer to \cite{Aizn-PAMS-2000,Alkhaleefah-PAMS-2019,Bhowmik-2018,Evdoridis-IM-2018,Kay & Pon & AASFM & 2019,Kayumov-CRACAD-2020,Liu-Ponnusamy-BMMS-2019} and the references therein.\vspace{1.5mm} 

\par Let $ \mathcal{H} $ be the class of all complex-valued harmonic functions $ f=h+\bar{g} $ defined on the unit disk $ \mathbb{D} $, where $ h $ and $ g $ are analytic $ \mathbb{D} $ with the normalization $ h(0)=h^{\prime}(0)-1=0 $ and $ g(0)=0 $. Let $ \mathcal{H}_0 $ be defined by $ 	\mathcal{H}_0=\{f=h+\bar{g}\in\mathcal{H} : g^{\prime}(0)=0\}. $ Then each $ f=h+\bar{g}\in \mathcal{H}_0 $ has the following form
\begin{equation}\label{e-1.1}
	h(z)=z+\sum_{n=2}^{\infty}a_nz^n\quad\mbox{and}\quad g(z)=\sum_{n=1}^{\infty}b_nz^n.
\end{equation}\vspace{1.5mm}

\par In $ 1977 $, Chichra \cite{Chichra-PAMS-1977} first introdcued the class $ \mathcal{W}(\alpha) $ consisting normalized analytic functions $ h $, satisfying the condition $ {\rm Re} \left(h^{\prime}(z)+\alpha zh^{\prime\prime}(z)>0\right) $ for $ z\in\mathbb{D} $ and $ \alpha\geq 0 $. Moreover, Chichra \cite{Chichra-PAMS-1977} has shown that functions in the class $ \mathcal{W}(\alpha) $ constitute a subclass of close-to-convex functions in $ \mathbb{D} $. In $ 2014 $, Nagpal and Ravichandran \cite{Nagpal-Ravinchandran-2014-JKMS} studied the following class 
\begin{equation*}
	\mathcal{W}^{0}_{\mathcal{H}}=\{f=h+\bar{g}\in \mathcal{H} :  {\rm Re}\left(h^{\prime}(z)+zh^{\prime\prime}(z)\right) > |g^{\prime}(z)+zg^{\prime\prime}(z)|\;\; \mbox{for}\; z\in\mathbb{D}\}
\end{equation*}
and obtained the coefficient bounds for the functions in the class $ \mathcal{W}^{0}_{\mathcal{H}} $. In $ 2019 $, Ghosh and Vasudevarao studied the class $ \mathcal{W}^{0}_{\mathcal{H}}(\alpha) $, where 
\begin{equation*}
	\mathcal{W}^{0}_{\mathcal{H}}(\alpha)=\{f=h+\bar{g}\in \mathcal{H} :  {\rm Re}\left(h^{\prime}(z)+\alpha zh^{\prime\prime}(z)\right) > |g^{\prime}(z)+\alpha zg^{\prime\prime}(z)|\;\; \mbox{for}\; z\in\mathbb{D}\}.
\end{equation*}
From the following results, it is easy to see that functions in the class $ \mathcal{W}^{0}_{\mathcal{H}}(\alpha) $ are univalent for $ \alpha\geq 0 $, and they are closely related to functions in $ \mathcal{W}(\alpha) $.
\begin{lem}\cite{Nirupam-MonatsMath-2019}\label{lem-1.2}
	The harmonic mapping $ f=h+\bar{g} $ belongs to $ \mathcal{W}^{0}_{\mathcal{H}}(\alpha) $ if, and only if, the analytic function $ F=h+\epsilon g $ belongs to $ \mathcal{W}(\alpha) $ for each $ |\epsilon|=1. $ 
\end{lem}
The coefficient bounds and the sharp growth estimates for functions in the class $ \mathcal{W}^{0}_{\mathcal{H}}(\alpha) $ have been studied in \cite{Nirupam-MonatsMath-2019}.
\begin{lem}\cite{Nirupam-MonatsMath-2019}\label{lem-1.3}
	Let $ f\in \mathcal{W}^{0}_{\mathcal{H}}(\alpha) $ for $ \alpha\geq 0 $ and be of the form \eqref{e-1.1}. Then for any $ n\geq 2 $,
	\begin{enumerate}
		\item[(i)] $ |a_n|+|b_n|\leq \displaystyle\frac{2}{\alpha n^2+(1-\alpha)n} $;\vspace{1.5mm}
		\item[(ii)] $  ||a_n|-|b_n||\leq \displaystyle\frac{2}{\alpha n^2+(1-\alpha)n} $; \vspace{1.5mm}
		\item[(iii)] $ |a_n|\leq \displaystyle\frac{2}{\alpha n^2+(1-\alpha)n} $.
	\end{enumerate}
All these inequalities are sharp for the function $ f=f_{\alpha} $ given by 
\begin{equation}\label{e-1.4}
	f_{\alpha}(z)=z+\sum_{n=2}^{\infty}\frac{2z^n}{\alpha n^2+(1-\alpha)n}.
\end{equation}
\end{lem}
\begin{lem}\cite{Nirupam-MonatsMath-2019}\label{lem-1.5}
	Let $ f\in \mathcal{W}^{0}_{\mathcal{H}}(\alpha) $ and be of the form \eqref{e-1.1} with $ 0<\alpha\leq 1 $. Then
	\begin{equation}\label{e-1.6}
		|z|+\sum_{n=2}^{\infty}\frac{2(-1)^{n-1}|z|^n}{\alpha n^2+(1-\alpha)n}\leq |f(z)|\leq |z|+\sum_{n=2}^{\infty}\frac{2|z|^n}{\alpha n^2+(1-\alpha)n}.
	\end{equation}
Both the inequalities are sharp for the function  $ f=f_{\alpha} $ given by \eqref{e-1.4}.
\end{lem}
\begin{lem}\cite{Nirupam-MonatsMath-2019}\label{lem-1.7}
		Let $ f\in \mathcal{W}^{0}_{\mathcal{H}}(\alpha) $ for $ \alpha\geq 0 $ and be of the form \eqref{e-1.1}. Then for $ n\geq 2 $, 
		\begin{equation*}
			|b_n|\leq\frac{1}{\alpha n^2+(1-\alpha)n}.
		\end{equation*} 
	The inequality is sharp for the function $ f=f^{*}_{\alpha} $ given by 
	\begin{equation}\label{e-1.8}
		f^{*}_{\alpha}(z) =z+\sum_{n=2}^{\infty}\frac{\overline{z^n}}{\alpha n^2+(1-\alpha)n}.
	\end{equation}
\end{lem}
In $ 2020 $, Vasudevarao and Halder \cite{Himadri-Vasu-P1} proved the following sharp Bohr radius for the class $ \mathcal{W}^{0}_{\mathcal{H}}(\alpha) $.
\begin{thm}\cite{Himadri-Vasu-P1}\label{th-11.99}
		Let $ f\in \mathcal{W}^{0}_{\mathcal{H}}(\alpha) $ for $ \alpha\geq 0 $ be of the form \eqref{e-1.1}. Then 
		\begin{equation*}
			|z|+\sum_{n=2}^{\infty}\left(|a_n|+|b_n|\right)|z|^n\leq d(f(0),\partial f(\mathbb{D}))
		\end{equation*} 
	for $ |z|=r\leq r_f $, where $ r_f $ is the unique root of 
	\begin{equation*}
		r+\sum_{n=2}^{\infty}\frac{2r^n}{\alpha n^2+(1-\alpha)n}=1+\sum_{n=2}^{\infty}\frac{2(-1)^{n-1}}{\alpha n^2+(1-\alpha)n}
	\end{equation*}
in $ (0,1) $. The radius $ r_f $ is the best possible.
\end{thm}
In $ 2020 $, Kayumov and Ponnusamy \cite{Kayumov-CRACAD-2020} obtained the following improved versions of Bohr inequality for analytic functions.
\begin{thm}\cite{Kayumov-CRACAD-2020}\label{th-1.9}
	Let $ f(z)=\sum_{n=0}^{\infty}a_nz^n $ be analytic in $ \mathbb{D} $, $ |f(z)|\leq 1 $ and $ S_r $ denotes the image of the subdisk $ |z|<r $ under the mapping $ f $. Then 
	\begin{equation*}
		B_1(r):=\sum_{n=0}^{\infty}|a_n|r^n+\frac{16}{9}\left(\frac{S_r}{\pi}\right)\leq 1\quad\mbox{for}\quad r\leq\frac{1}{3},
	\end{equation*} 
and the number $ 1/3 $ and $ 16/9 $ cannot be improved. Moreover, 
	\begin{equation*}
	B_2(r):=|a_0|^2+\sum_{n=1}^{\infty}|a_n|r^n+\frac{8}{9}\left(\frac{S_r}{\pi}\right)\leq 1\quad\mbox{for}\quad r\leq\frac{1}{2},
\end{equation*} 
and the number $ 1/2 $ and $ 8/9 $ cannot be improved.
\end{thm}
In analogy to Bohr radius, Bohr-Rogosinski radius has also been defined (see \cite{Rogosinski-1923}) which is describes as follows: If $ f\in\mathcal{B} $, then for $ N\geq 1 $, we have $ |S_N(z)|<1 $ in the disk $ \mathbb{D}_{1/2} $ and this radius is sharp, where $ S_N(z)=\sum_{n=0}^{N}a_nz^n $ denotes the partial sum of $ f $. There is a relavent quantity, which we call the Bohr-Rogosinski sum $ R^f_N(z) $ of $ f $ defined by 
\begin{equation}\label{e-11.12}
	R^f_N(z):=|f(z)|+\sum_{n=N}^{\infty}|a_n|r^n,\;\; |z|=r.
\end{equation}  
\par It is important to note that for $ N=1$, the quantity \eqref{e-11.12} is reduces to the classical Bohr sum in which $ f(0) $ is replaced by $ |f(z)| $. In $ 2017 $, Kayumov and Ponnusamy \cite{Kayumov-2017} proved the following interesting result on Bohr-Rogosinski radius for the analytic functions.
\begin{thm}\cite{Kayumov-2017}\label{th-11.1111}
	Suppose that $ f(z)=\sum_{n=0}^{\infty}a_nz^n $ is analytic in the unit disk $ \mathbb{D} $ and $ |f(z)|<1 $ in $ \mathbb{D} $. Then 
	\begin{equation*}
		|f(z)|+\sum_{n=N}^{\infty}|a_n|r^n\leq 1\;\;\mbox{for}\;\; r\leq R_N,
	\end{equation*}
where $ R_N $ is the positive root of the equation $ 2(1+r)r^N-(1-r^2)=0 $. The radius $ R_N $ is the best possible. Moreover, 
\begin{equation*}
	|f(z)|^2+\sum_{n=N}^{\infty}|a_n|r^n\leq 1;\;\mbox{for}\;\; r\leq R^{\prime}_N,
\end{equation*}
where $ R^{\prime}_N $ is the positive root of the equation $ (1+r)r^N-(1-r^2)=0 $. The radius $ R^{\prime}_N $ is the best possible.
\end{thm}
The main aim of this paper is to study the Bohr phenomenon for the class $ \mathcal{W}^0_{\mathcal{H}}(\alpha) $ for $ \alpha\geq 0 $. We find the improved Bohr radius, refined Bohr radius as well as Bohr-Rogosinski inequality for functions in the class $ \mathcal{W}^0_{\mathcal{H}}(\alpha) $. In Section 2, we state the main results of the paper, and in Section 3, we prove all the main results.
\section{Main results}
 Before stating the main results of this paper, we recall the definiton of dilogarithm. The dilogarithm function $ {\rm Li}_2(z) $ is defined by the power series
\begin{equation*}
	{\rm Li}_2(z)=\sum_{n=1}^{\infty}\frac{z^n}{n^2}\;\; \mbox{for}\;\; |z|<1.
\end{equation*}
The definition and the name, of course, come from the analogy with the Taylor series of the ordinary logarithm around $ 1 $, 
\begin{equation*}
	-\log(1-z)=\sum_{n=1}^{\infty}\frac{z^n}{n}\;\;\mbox{for}\;\; |z|<1,
\end{equation*}
which leads similarly to the defintion of the polylogarithm
\begin{equation*}
	{\rm Li}_m(z)=\sum_{n=1}^{\infty}\frac{z^n}{n^m}\;\; \mbox{for}\;\; |z|<1,\;\; m=1, 2, 3, \ldots.
\end{equation*}
\begin{figure}[!htb]
	\begin{center}
		\includegraphics[width=0.47\linewidth]{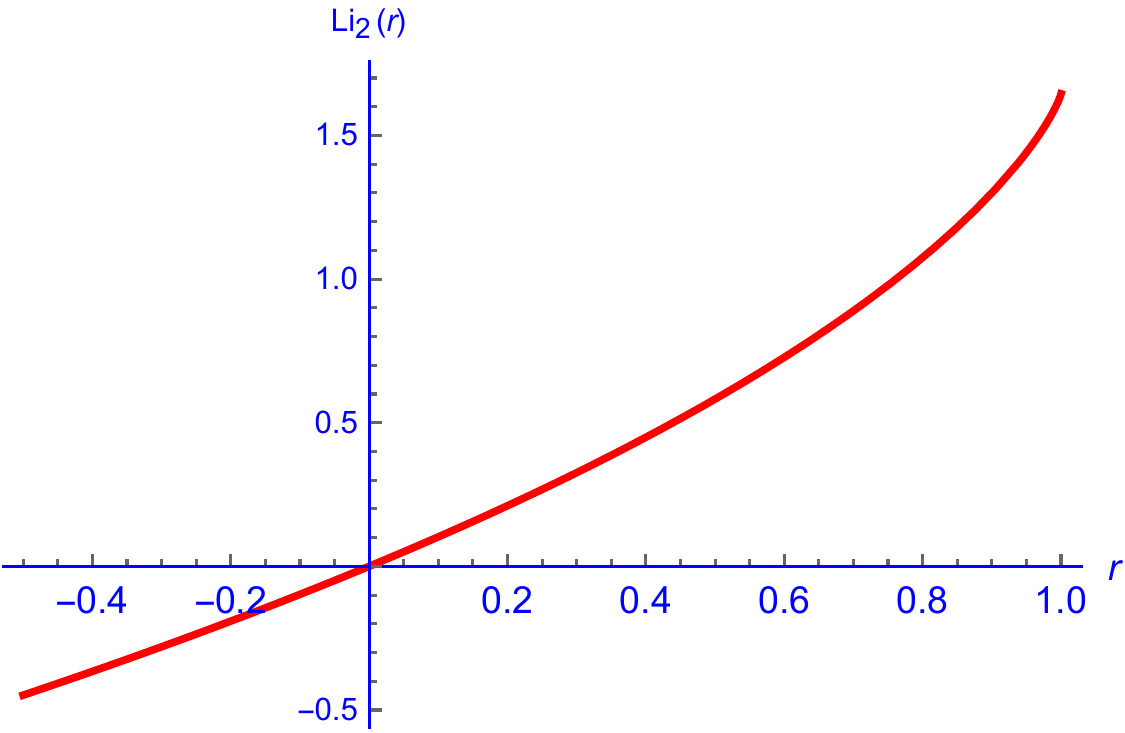}
	\end{center}
	\caption{Graph of $ {\rm Li}_2(r) $ for $ r\in (-0.5,1) $.}
\end{figure}
In addition, the relation 
\begin{equation*}
	\frac{d}{dz}\left({\rm Li}_m(z)\right)=\frac{1}{z}{\rm Li}_{m-1}(z)\;\; \mbox{for}\;\; m\geq 2
\end{equation*}
is obvious and leads by induction to the extension of the domain of definition of $ {\rm Li}_m $ to the cut plane $ \mathbb{C}\setminus [1, \infty) $. In particular, the analytic continuation of the dilogarithm is given by
\begin{equation*}
	{\rm Li}_2(z)=-\int_{0}^{z}\log (1-u)\frac{du}{u}\;\; \mbox{for}\;\; z\in\mathbb{C}\setminus [1,\infty).
\end{equation*}
\par Using Lemma \ref{lem-1.3} and Lemma \ref{lem-1.5}, considering power of the coefficient $ |a_n|+|b_n| $, we prove the following sharp Bohr radius for the class $ \mathcal{W}^{0}_{\mathcal{H}}(\alpha) $.
\begin{thm}\label{th-2.1}
	Let $ f\in \mathcal{W}^{0}_{\mathcal{H}}(\alpha) $ for $ \alpha\geq 0 $ be of the form \eqref{e-1.1}. Then for $ n\geq 2 $ and $ p\geq 1 $,
	\begin{equation}\label{e-2.2}
		|z|+\sum_{n=2}^{\infty}(|a_n|+|b_n|)|z|^n+\sum_{n=2}^{\infty}(|a_n|+|b_n|)^p|z|^{pn}\leq d(f(0),\partial f(\mathbb{D}))
	\end{equation}
for $ |z|=r\leq r_p(\alpha) $, where $ r_p(\alpha) $ is the unique root of 
\begin{equation}\label{e-2.3}
	r+\sum_{n=2}^{\infty}\frac{2r^n}{\alpha n^2+(1-\alpha)n}+\sum_{n=2}^{\infty}\left(\frac{2r^n}{\alpha n^2+(1-\alpha)n}\right)^p=1+\sum_{n=2}^{\infty}\frac{2(-1)^{n-1}}{\alpha n^2+(1-\alpha)n}.
\end{equation}
Here $ r_p(\alpha) $ is the best possible.
\end{thm}
For a particular choice of $ \alpha=1/2 $ and $ p=2 $, we have the following corollary of Theorem \ref{th-2.1}.
\begin{cor}\label{cor-2.4}
	Let $ f\in \mathcal{W}^{0}_{\mathcal{H}}(1/2) $ and be of the form \eqref{e-1.1}. Then for $ |z|=r $,\begin{equation*}
		|z|+\sum_{n=2}^{\infty}(|a_n|+|b_n|)|z|^n+\sum_{n=2}^{\infty}(|a_n|+|b_n|)^2|z|^{2n}\leq d(f(0), \partial f(\mathbb{D}))
	\end{equation*}
	for $ r_2\leq r_2(0.5)\approx 0.399085 $, where $ r_2(0.5) $ is the unique root of the equation
\begin{align*}
		&-4r^2-r-39+\frac{16}{r^2}\left((1+r^2){\rm Li_2(r^2)}-32(1-r^2)\log(1-r^2)\right)\\&\quad\quad+\frac{4}{r}\left((1-r)\log(1-r)\right)-8\log 2=0
\end{align*}
	in $ (0,1) $. Here $ r_2(0.5)\approx 0.399085 $ is the best possible.
\end{cor}
\begin{cor}\label{cor-2.5}
		Let $ f\in \mathcal{W}^{0}_{\mathcal{H}}(1/2) $ be of the form \eqref{e-1.1}. Then for $ |z|=r $,\begin{equation*}
		|z|+\sum_{n=2}^{\infty}(|a_n|+|b_n|)^2|z|^{2n}\leq d(f(0), \partial f(\mathbb{D}))
	\end{equation*}
	for $ r\leq r^*(0.5)\approx 0.512331 $, where $ r^*(0.5) $ is the unique root of the equation
	\begin{align*}
		&-4r^2+r+\frac{16(1+r^2){\rm Li_2}(r^2)}{r^2}-\frac{32(1-r^2)\log(1-r^2)}{r^2}-43-8\log 2=0
	\end{align*}
	in $ (0,1) $. Here $ r^*(0.5)\approx 0.512331 $ is the best possible.
\end{cor}
\begin{figure}[!htb]
	\begin{center}
		\includegraphics[width=0.47\linewidth]{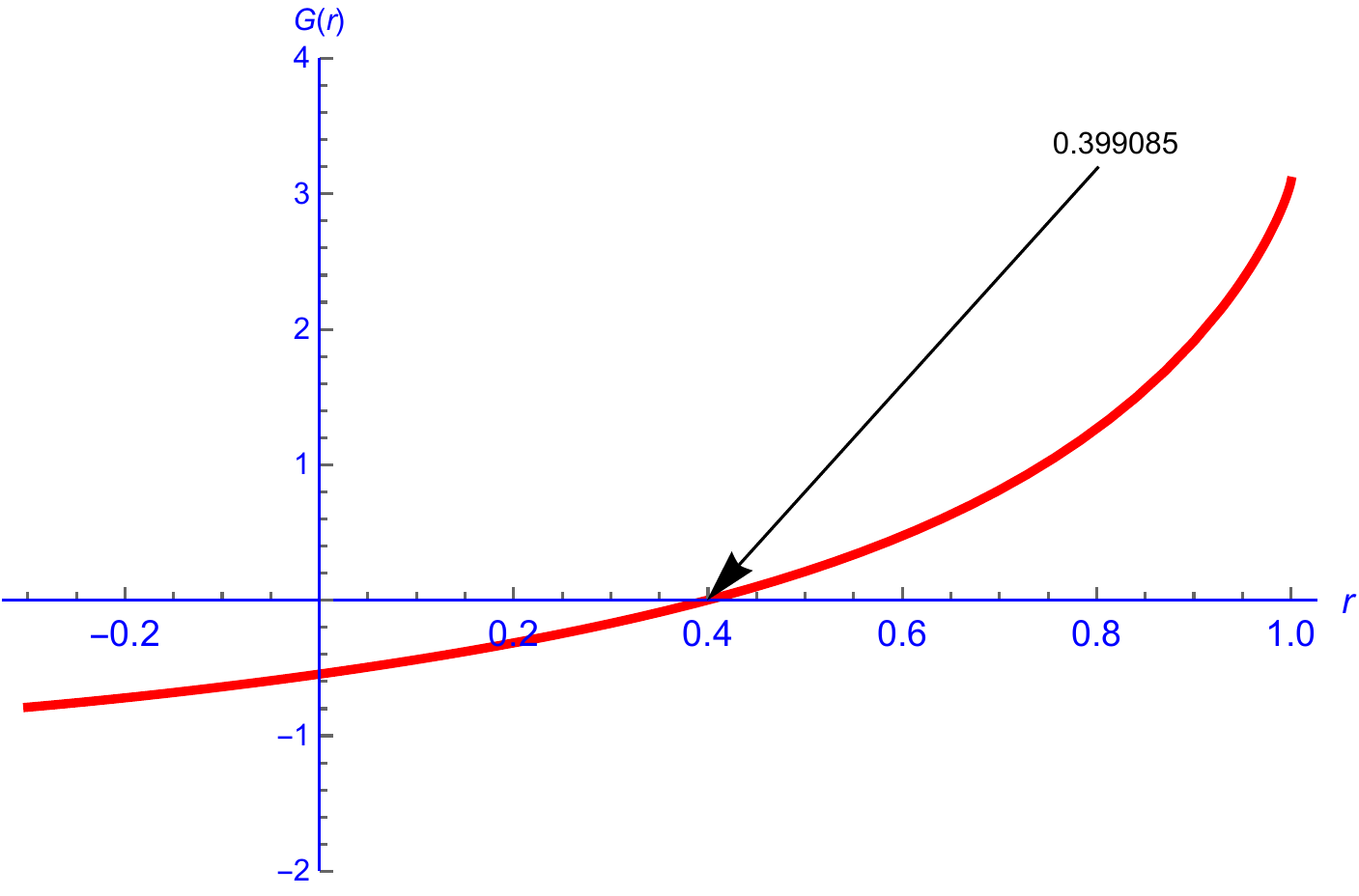}\;\;  \includegraphics[width=0.47\linewidth]{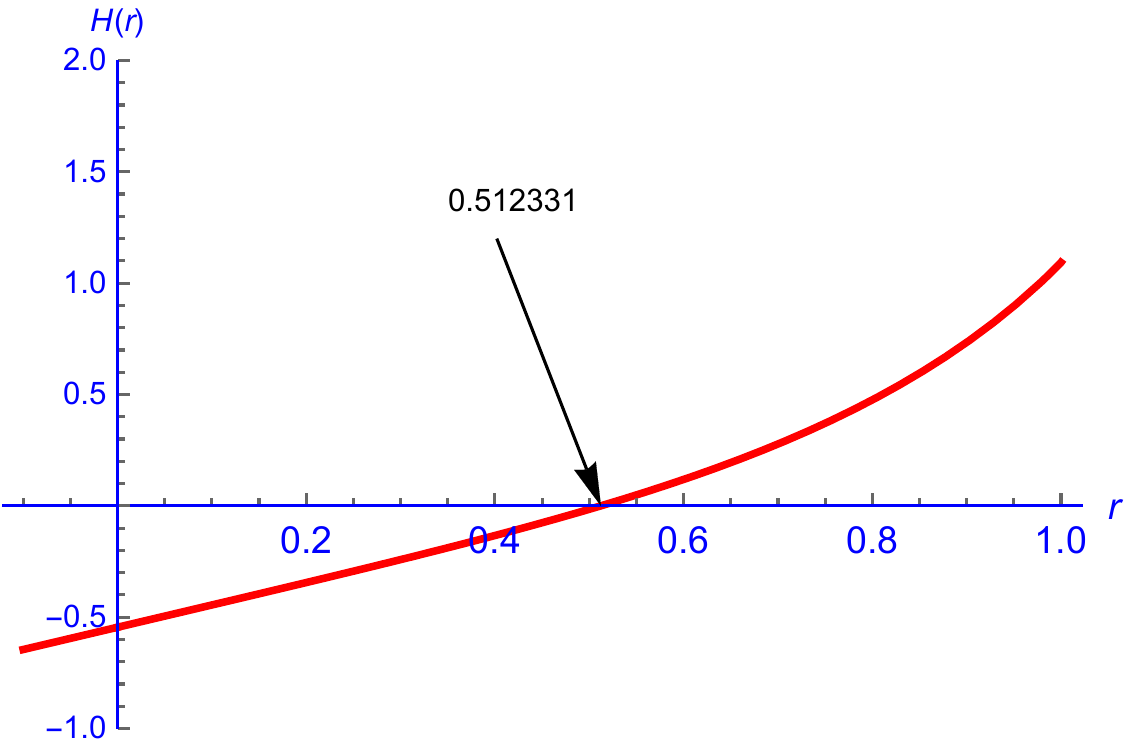}
	\end{center}
	\caption{The radii $  r_2(0.5)\approx 0.399085 $ and $ r^*(0.5)\approx 0.512331 $ in Corollary \ref{cor-2.4} and Corollary \ref{cor-2.5} respectively.}
\end{figure}
\begin{rem}
	It is easy to see that Theorem \ref{th-11.99} is a particular case of Theorem \ref{th-2.1}. In fact, Theorem \ref{th-2.1} is an improved version of Theorem \ref{th-11.99}.
\end{rem}
\begin{rem}
	For particular values of $ \alpha $ and $ p $ in Theorem \ref{th-2.1}, we discuss the Bohr radii for the class $ \mathcal{W}^{0}_{\mathcal{H}}(\alpha) $.
	\begin{enumerate}
		\item[(i)] For $ \alpha=0 $, \eqref{e-2.3} reduces to the following
		\begin{equation}\label{e-2.4}
			r-2(r+\log(1-r))-2^p(r^p+\log(1-r^p))+1-2\log 2=0.
		\end{equation}
	A simple computation shows that $ r_1(0)\approx 0.243755  $, $ r_3(0)\approx 0.284007  $, $ r_8(0)\approx 0.285194  $, $ r_{17}(0)\approx 0.285194  $. In fact $ r_p(0)\approx 0.285194  $ for $ p\geq 8 $.\vspace{1.5mm}
	\item[(ii)] For $ \alpha=1/2 $, the radius are $ r_1(0.5)\approx 0.347966 $, $ r_2(0.5)\approx 0.399085 $.\vspace{1.5mm}
	\item[(iii)] For $ \alpha=1 $, \eqref{e-2.3} takes the form
\begin{equation}\label{e-2.5}
	r+2(-r+{\rm Li_2}(r))+2^p(-r^p+{\rm Li_{2p}}(r^p))+1-\frac{\pi^2}{6}=0.
\end{equation}
A simple computation shows that $ r_1(1)\approx 0.422211  $, $ r_2(1)\approx 0.480812  $, $ r_3(1)\approx 0.487911  $, $ r_{4}(1)\approx 0.488711  $, $ r_{5}(1)\approx 0.488874  $, $ r_{6}(1)\approx 0.488886  $. In fact, $ r_p(1)\approx 0.488888  $ for $ p\geq 7 $.
	\end{enumerate}
\end{rem}
 \begin{figure}[!htb]
	\begin{center}
		\includegraphics[width=0.45\linewidth]{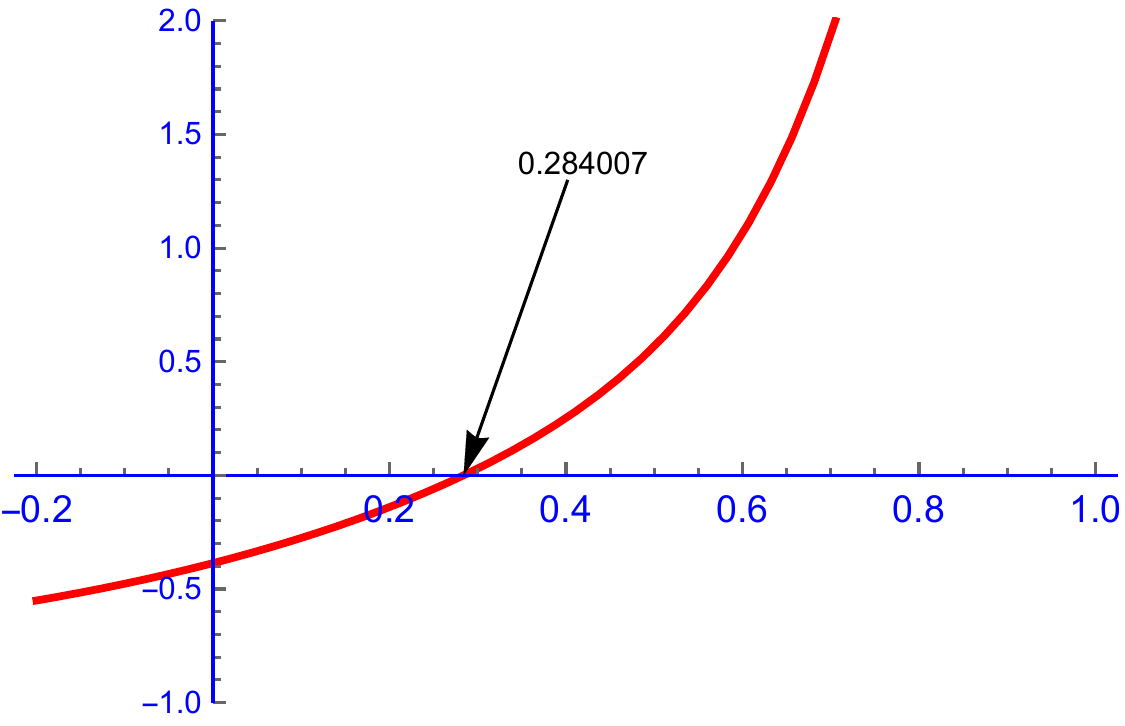}\;\;  \includegraphics[width=0.45\linewidth]{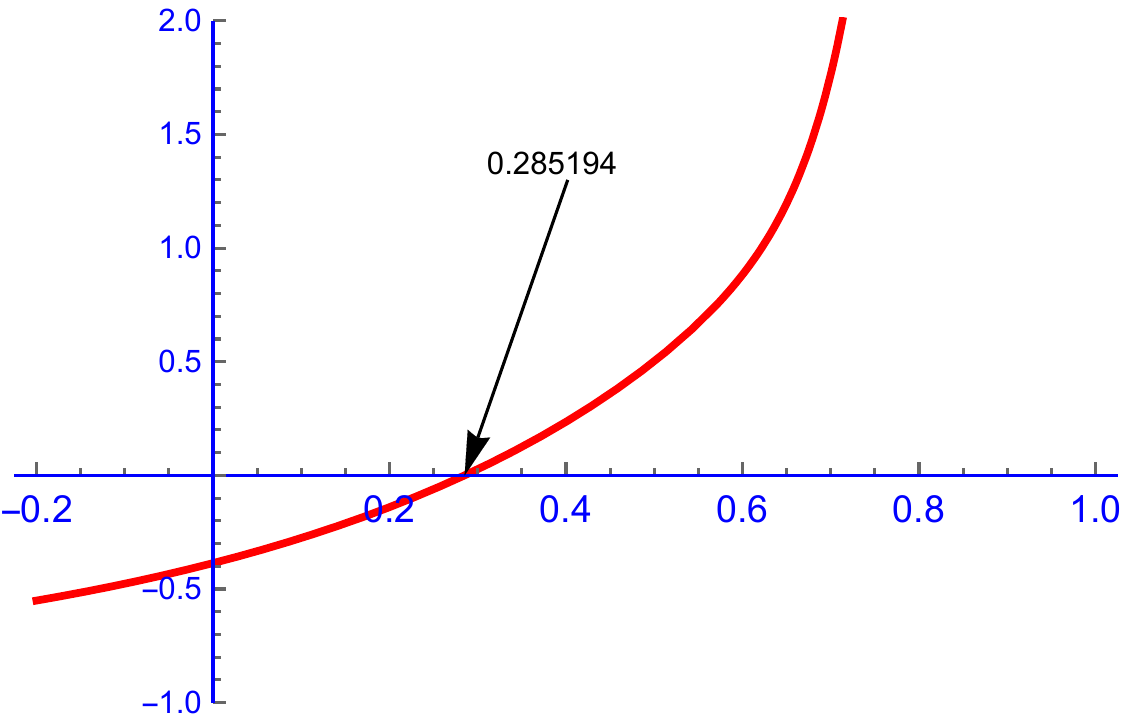}
	\end{center}
	\caption{The radii $  r_3(0)\approx 0.284007 $ and $ r_8(0)\approx 0.285194 $ which are roots of the equation \eqref{e-2.4}.}
\end{figure}
\begin{figure}[!htb]
	\begin{center}
		\includegraphics[width=0.45\linewidth]{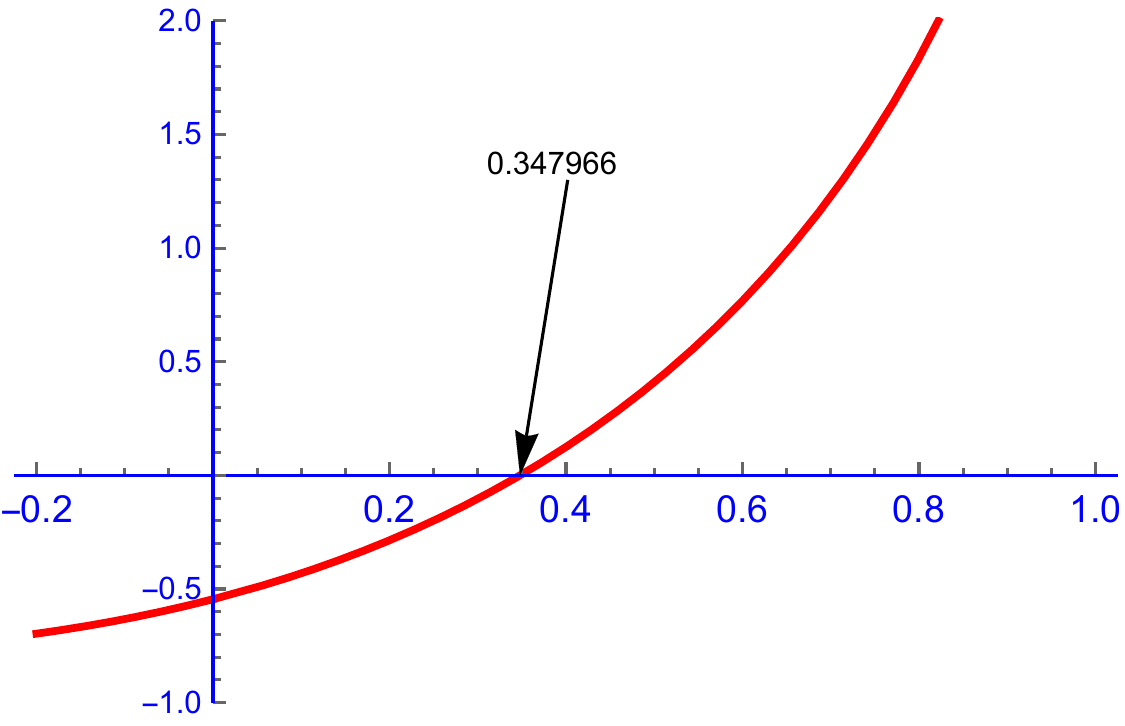}\;\;  \includegraphics[width=0.45\linewidth]{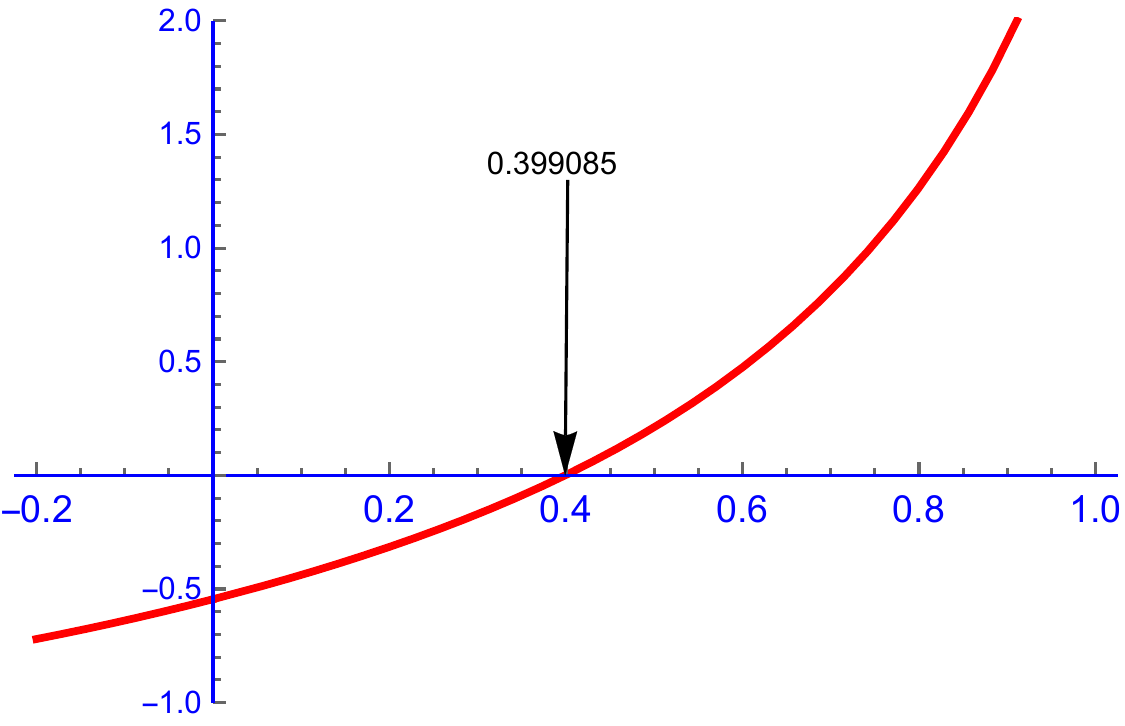}
	\end{center}
	\caption{The radii $  r_1(0.5)\approx 0.347966 $ and $ r_2(0.5)\approx 0.399085 $ for $ \alpha=1/2 $.}
\end{figure}
\begin{figure}[!htb]
	\begin{center}
		\includegraphics[width=0.45\linewidth]{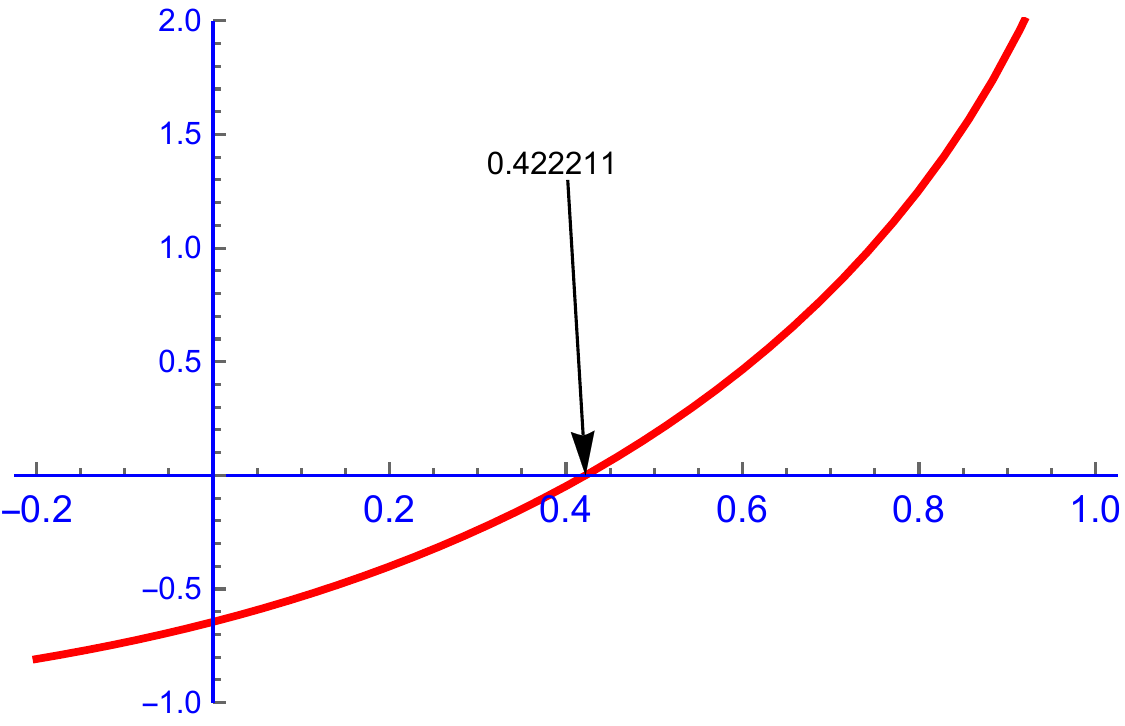}\;\;  \includegraphics[width=0.45\linewidth]{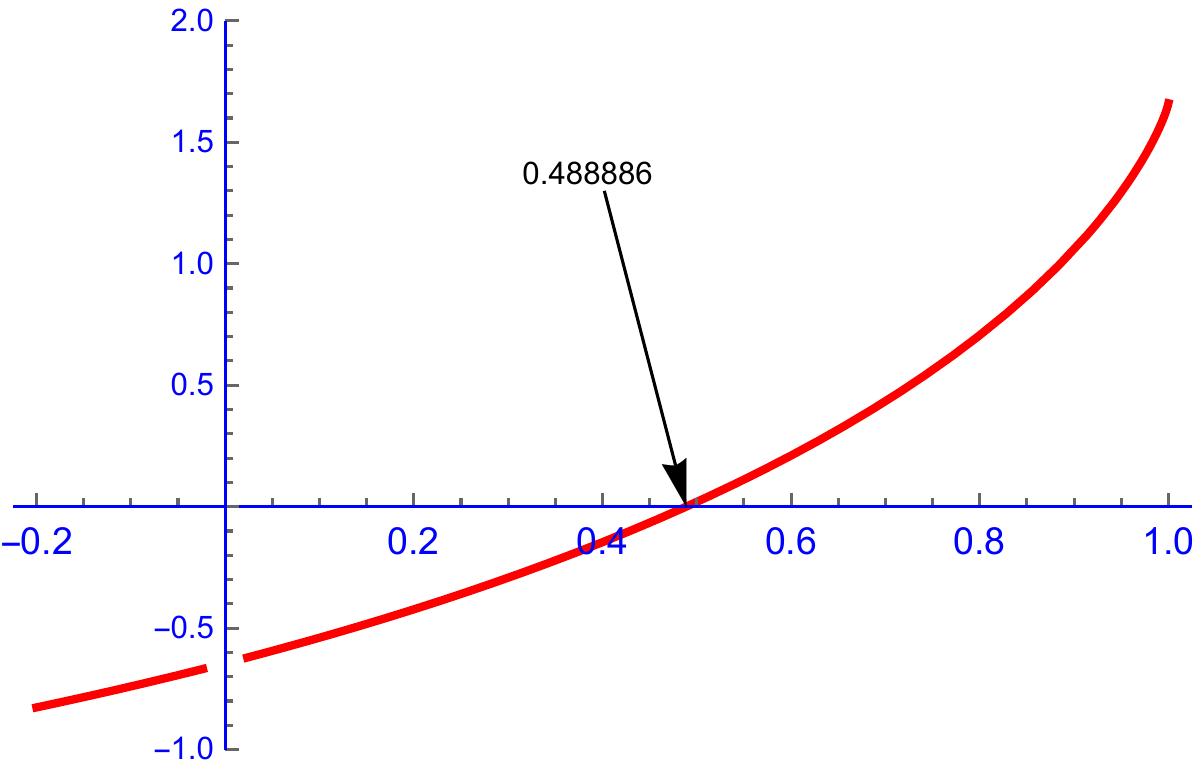}
	\end{center}
	\caption{The radii $  r_1(1)\approx 0.422211 $ and $ r_6(0)\approx 0.488886 $ which are roots of the equation \eqref{e-2.5}.}
\end{figure}

\begin{thm}\label{th-2.6}
	Let $ f\in \mathcal{W}^{0}_{\mathcal{H}}(\alpha) $ for $ \alpha\geq 0 $ be of the form \eqref{e-1.1}. Then for $ n\geq 2 $,
	\begin{equation}
		|f(z)|+\sum_{n=2}^{\infty}(|a_n|+|b_n|)|z|^n\leq d(f(0),\partial f(\mathbb{D}))
	\end{equation}
	for $ |z|=r\leq r_*(\alpha) $, where $ r_*(\alpha) $ is the unique root of 
	\begin{equation*}
		r+\sum_{n=2}^{\infty}\frac{4r^n}{\alpha n^2+(1-\alpha)n}-1-\sum_{n=2}^{\infty}\frac{2(-1)^{n-1}}{\alpha n^2+(1-\alpha)n}=0.
	\end{equation*}
	Here $ r_*(\alpha) $ is the best possible. 
\end{thm}

For a particular choice of $ \alpha $, we have the following corolarry of Theorem \ref{th-2.6}.
\begin{cor}\label{cor-2.10}
	Let $ f\in \mathcal{W}^{0}_{\mathcal{H}}(1/2) $ be of the form \eqref{e-1.1}. Then for $ |z|=r $,
	\begin{equation*}
		|f(z)|+\sum_{n=2}^{\infty}(|a_n|+|b_n|)|z|^n\leq d(f(0), \partial f(\mathbb{D}))
	\end{equation*}
	for $ r\leq r_*(0.5)\approx 0.347966 $, where $ r_*(0.5) $ is the unique root of the equation
	\begin{equation*}
		\frac{8}{r}(1-r)\log(1-r)-3r+13-8\log 2=0
	\end{equation*}
	in $ (0,1) $. Here $ r_*(0.5)\approx 0.347966 $ is the best possible.
\end{cor}
\begin{rem}
	For $ p=1 $, the Bohr radius in Theorem \ref{th-2.1} coincides with the Bohr radius in Theorem \ref{th-2.6}. In addition, for $ \alpha=1/2 $ and $ p=1 $, Bohr radius in Theorem \ref{th-2.1} coincides with that of in Corollary \ref{cor-2.10}.
\end{rem}
For the harmonic analogue of Theorem \ref{th-11.1111}, we prove the following sharp Bohr-Rogosinski radius for the class $ f\in \mathcal{W}^{0}_{\mathcal{H}}(\alpha) $.
\begin{thm}\label{th-2.11}
	Let $ f\in \mathcal{W}^{0}_{\mathcal{H}}(\alpha) $ and be of the form \eqref{e-1.1}. Then for integers $ m\geq 1 $ $ n, N\geq 2 $, we have
\begin{equation}
	|f(z^m)|+\sum_{n=N}^{\infty}(|a_n|+|b_n|)|z|^n\leq d(f(0),\partial f(\mathbb{D}))
\end{equation}
for $ |z|=r\leq R_{m,N}(\alpha) $, where $ R_{m,N}(\alpha) $ is the unique root of 
\begin{equation*}
	r^m+\sum_{n=2}^{\infty}\frac{2r^{mn}}{\alpha n^2+(1-\alpha)n}+\sum_{n=N}^{\infty}\frac{2r^{n}}{\alpha n^2+(1-\alpha)n}-1-\sum_{n=2}^{\infty}\frac{2(-1)^{n-1}}{\alpha n^2+(1-\alpha)n}=0.
\end{equation*}
Here $ R_{m,N}(\alpha) $ is the best possible.
\end{thm}
Our primary objective is now to generalize the harmonic versions of Theorem \ref{th-1.9} for functions in the class $ \mathcal{W}^{0}_{\mathcal{H}}(\alpha) $.
\begin{thm}\label{th-2.13}
	Let $ f\in \mathcal{W}^{0}_{\mathcal{H}}(\alpha)  $ be given by \eqref{e-1.1}. Then 
	\begin{enumerate}
		\item[(i)]  
		\begin{equation*}
			|z|+\sum_{n=2}^{\infty}\left(|a_n|+|b_n|\right)|z|^n+\frac{S_r}{\pi}\leq d(f(0),\partial f(\mathbb{D}))
		\end{equation*}
	for $ |z|=r\leq r_f(\alpha) $, where $ r_f(\alpha) $ is the unique root of the equation
\begin{align*}
		&r^2+r+\sum_{n=2}^{\infty}\frac{2r^n}{\alpha n^2+(1-\alpha)n}+\sum_{n=2}^{\infty}\frac{4nr^{2n}}{(\alpha n^2+(1-\alpha)n)^2}-1\\&\quad\quad-\sum_{n=2}^{\infty}\frac{2(-1)^{n-1}}{\alpha n^2+(1-\alpha)n}=0
\end{align*}
in $ (0,1) $. Here $ r_f(\alpha) $ is the best possible.\vspace{1.5mm}
\item[(ii)] 
\begin{equation*}
	|f(z)|^2+\sum_{n=2}^{\infty}\left(|a_n|+|b_n|\right)|z|^n+\left(\frac{S_r}{\pi}\right)^2\leq d(f(0),\partial f(\mathbb{D}))
\end{equation*}
for $ |z|=r\leq r^{*}_f(\alpha) $, where $ r^{*}_f(\alpha) $ is the unique root of the equation
\begin{align*}
	&\left(r+\sum_{n=2}^{\infty}\frac{2r^n}{\alpha n^2+(1-\alpha)n}\right)^2+\sum_{n=2}^{\infty}\frac{2r^n}{\alpha n^2+(1-\alpha)n}+\left(r^2+\sum_{n=2}^{\infty}\frac{4nr^{2n}}{(\alpha n^2+(1-\alpha)n)^2}\right)^2\\&\quad\quad-1-\sum_{n=2}^{\infty}\frac{2(-1)^{n-1}}{\alpha n^2+(1-\alpha)n}=0
\end{align*}
in $ (0,1) $.  Here $ r^{*}_f(\alpha) $ is the best possible.
	\end{enumerate}
\end{thm} 
For a particular choice of $ \alpha $, we have the following corollary of Theorem \ref{th-2.13}.
\begin{cor}\label{cor-2.14}
		Let $ f\in \mathcal{W}^{0}_{\mathcal{H}}(\alpha)  $ be given by \eqref{e-1.1}. Then
	\begin{enumerate}
		\item[(i)]   for $ \alpha=1/2 $,
		\begin{equation*}
			|z|+\sum_{n=2}^{\infty}\left(|a_n|+|b_n|\right)|z|^n+\frac{S_r}{\pi}\leq d(f(0),\partial f(\mathbb{D}))
		\end{equation*}
		for $ |z|=r\leq r_f(\alpha)\approx 0.600881 $, where $ r_f(\alpha) $ is the unique root of the equation
		\begin{align*}
		&\frac{4}{r}(1-r)\log (1-r)+\frac{16}{r^2}(1-r^2)\log(1-r^2)-\frac{16}{r^2}{\rm Li}_2(r^2)\\&\quad\quad-3r^2-r+29+8\log 2=0
		\end{align*}
		in $ (0,1) $. Here $ r_f(\alpha)\approx 0.600881 $ is the best possible.\vspace{1.5mm}
		\item[(ii)]  for $ \alpha=1/2 $,
		\begin{equation*}
			|f(z)|+\sum_{n=2}^{\infty}\left(|a_n|+|b_n|\right)|z|^n+\frac{S_r}{\pi}\leq d(f(0),\partial f(\mathbb{D}))
		\end{equation*}
		for $ |z|=r\leq r^{*}_f(\alpha)\approx 0.302059 $, where $ r^{*}_f(\alpha) $ is the unique root of the equation
		\begin{align*}
			&\frac{8}{r}(1-r)\log (1-r)+\frac{16}{r^2}(1-r^2)\log (1-r^2)-\frac{16}{r^2}{\rm Li}_2(r^2)\\&\quad\quad -3r^2-3r+45-8\log 2=0
		\end{align*}
		in $ (0,1) $. Here $ r^{*}_f(\alpha)\approx 0.302059 $ is the best possible.
	\end{enumerate}
\end{cor}
\begin{figure}[!htb]
	\begin{center}
		\includegraphics[width=0.45\linewidth]{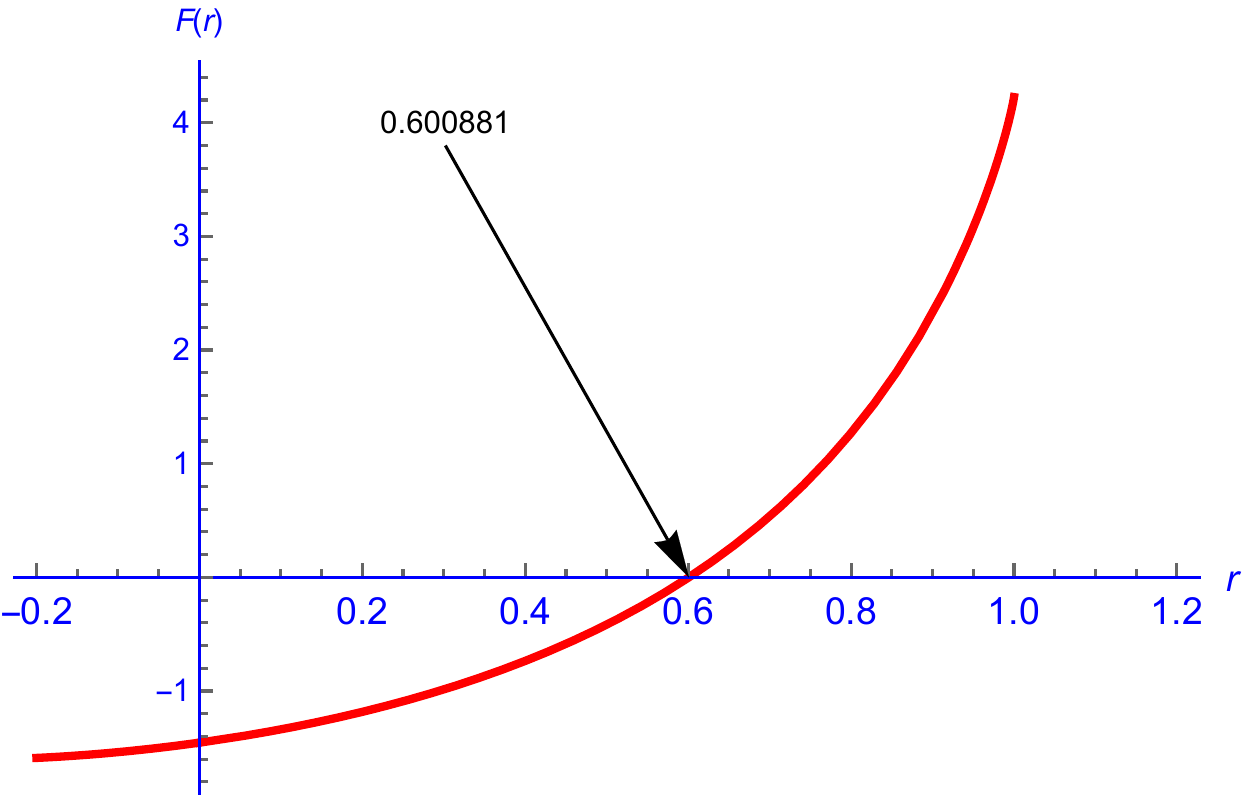}\;\;  \includegraphics[width=0.45\linewidth]{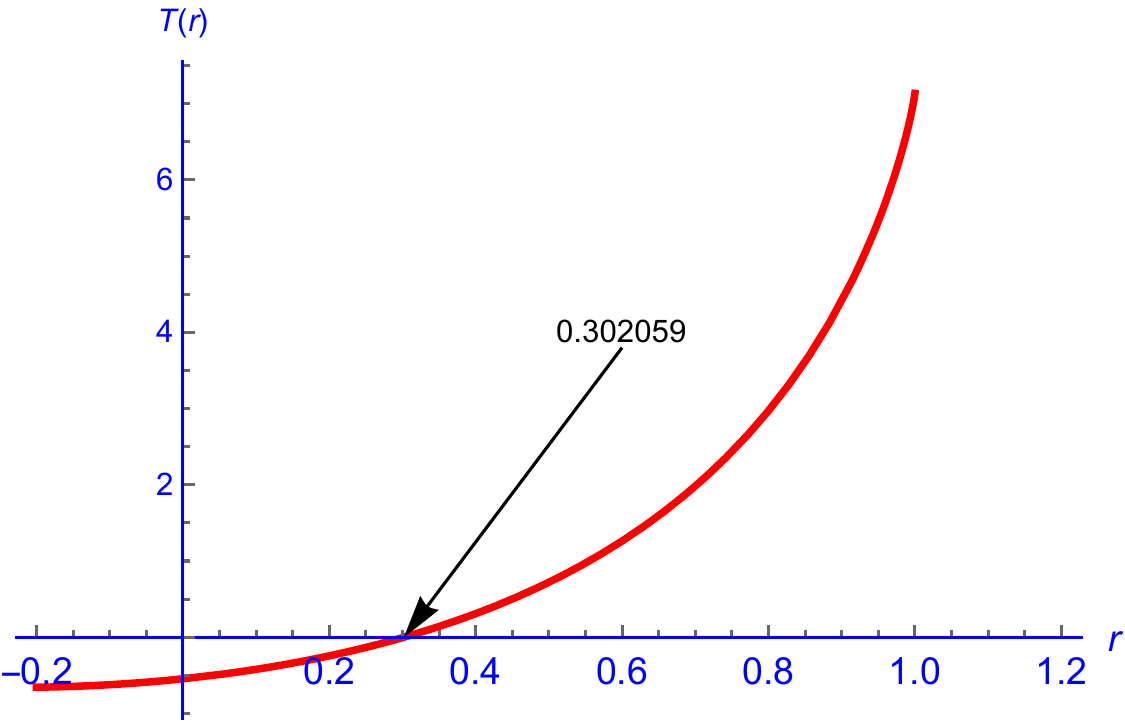}
	\end{center}
	\caption{The radii $  r_f(\alpha)\approx 0.600881 $ and $ r^*_f(\alpha)\approx 0.302059 $.}
\end{figure}
\begin{thm}\label{th-2.15}
	Let $ f\in \mathcal{W}^{0}_{\mathcal{H}}(\alpha)  $ be given by \eqref{e-1.1}. Then 
	\begin{enumerate}
		\item[(i)] 
		\begin{equation*}
			|z|+|h(z)|+\sum_{n=2}^{\infty}|a_n||z|^n\leq d(f(0),\partial f(\mathbb{D}))
		\end{equation*}
		for $ |z|=r\leq r_h(\alpha) $, where $ r_h(\alpha) $ is the unique root of the equation
		\begin{equation}\label{e-2.16}
			2r+\sum_{n=2}^{\infty}\frac{4r^n}{\alpha n^2+(1-\alpha)n}-1-\sum_{n=2}^{\infty}\frac{2(-1)^{n-1}}{\alpha n^2+(1-\alpha)n}=0
		\end{equation}
		in $ (0,1) $. Here $ r_h(\alpha) $ is the best possible.
	\end{enumerate}
\item[(ii)] 
\begin{equation*}
	|z|+|g(z)|+\sum_{n=2}^{\infty}|a_n||z|^n\leq d(f(0),\partial f(\mathbb{D}))
\end{equation*}
for $ |z|=r\leq r_g(\alpha) $, where $ r_g(\alpha) $ is the unique root of the equation
\begin{equation}\label{e-2.17}
	r+\sum_{n=2}^{\infty}\frac{2r^n}{\alpha n^2+(1-\alpha)n}-1-\sum_{n=2}^{\infty}\frac{2(-1)^{n-1}}{\alpha n^2+(1-\alpha)n}=0
\end{equation}
in $ (0,1) $. Here $ r_g(\alpha) $ is the best possible.
\end{thm}
\begin{rem}
	A simple computation shows that $ r_h(0)\approx 0.1632 $, $ r_h(1/2)\approx 0.2319 $, $ r_h(1/3)\approx 0.2119 $, $ r_h(1/4)\approx 0.2009 $, $ r_h(1/8)\approx 0.1829 $ and $ r_g(0)\approx 0.2852 $, $ r_g(1/2)\approx 0.4057 $, $ r_g(1/3)\approx 0.3704 $, $ r_g(1/4)\approx 0.351 $, $ r_g(1/8)\approx 0.3195 $.
\end{rem}
\begin{cor}\label{cor-2.18}
	Let $ f\in\mathcal{W}^0_{\mathcal{H}}(\alpha) $ be given by \eqref{e-1.1}.. Then for $ \alpha=1/2 $ and $ |z|=r $, 
	\begin{equation*}
		|g(z)|+\sum_{n=2}^{\infty}|a_n||z|^n\leq d(f(0),\partial f(\mathbb{D}))
	\end{equation*}
for $ |z|=r\leq r^*_g(1/2)\approx 0.794054 $. Here $ r^*_g(1/2) $ is the best possible.
\end{cor}
\begin{figure}[!htb]
	\begin{center}
		\includegraphics[width=0.50\linewidth]{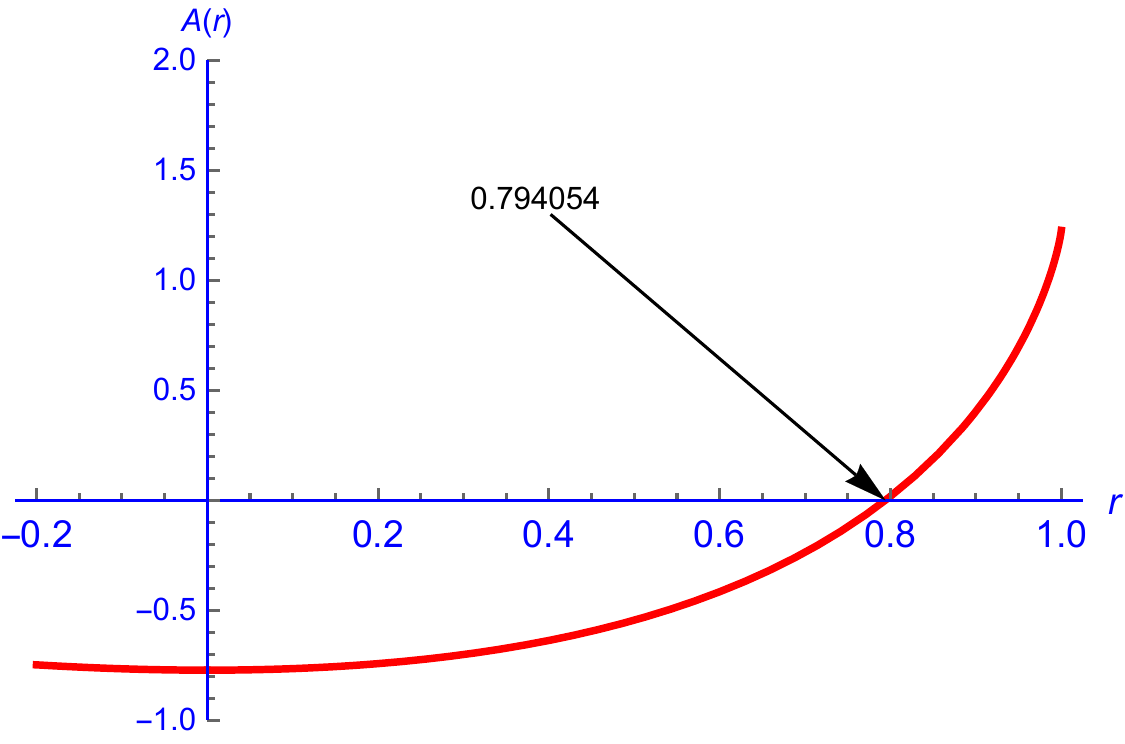}
	\end{center}
	\caption{The figure illustrates the radius $ r_g(1/2)\approx 0.704287 $.}
\end{figure}

\section{Proof of the main results}
\begin{proof}[\bf Proof of Theorem \ref{th-2.1}]
		Let $ f\in \mathcal{W}^{0}_{\mathcal{H}}(\alpha)  $ be given by \eqref{e-1.1}. Then in view of Lemma \ref{lem-1.3} and Lemma \ref{lem-1.5}, it is evident that the Euclidean distance $ d(f(0),\partial f(\mathbb{D})) $ between $ f(0) $ and the boundary of $ f(\mathbb{D}) $ is 
		\begin{equation}\label{e-3.1}
			d(f(0),\partial f(\mathbb{D}))=\liminf_{|z|\rightarrow 1}|f(z)-f(0)|\geq 1+\sum_{n=2}^{\infty}\frac{2(-1)^{n-1}}{\alpha n^2+(1-\alpha)n}.
		\end{equation} 
	Let $ J_1 : [0,1]\rightarrow\mathbb{R} $ be defined by 
	\begin{align}\label{ee-33.22}
		J_1(r)&=r+\sum_{n=2}^{\infty}\frac{2r^n}{\alpha n^2+(1-\alpha)n}+\sum_{n=2}^{\infty}\left(\frac{2r^n}{\alpha n^2+(1-\alpha)n}\right)^p-1\\&\nonumber\quad\quad -\sum_{n=2}^{\infty}\frac{2(-1)^{n-1}}{\alpha n^2+(1-\alpha)n}.
	\end{align}
Clearly, the function $ J_1(r) $ is continuous on $ [0,1] $ and differentiable on $ (0,1) $. Since, 
\begin{equation*}
	\bigg|\sum_{n=2}^{\infty}\frac{(-1)^{n-1}}{\alpha n^2+(1-\alpha)n}\bigg|\leq \frac{1}{2}\;\; \mbox{for}\;\; n\geq 2.
\end{equation*}
it is easy to see that
\begin{equation*}
	J_1(0)=-1-\sum_{n=2}^{\infty}\frac{2(-1)^{n-1}}{\alpha n^2+(1-\alpha)n}<0.
\end{equation*}
On the other hand, since  
\begin{equation*}
	\sum_{n=2}^{\infty}\frac{1}{\alpha  n^2+(1-\alpha)n}\geq \sum_{n=2}^{\infty}\frac{(-1)^{n-1}}{\alpha  n^2+(1-\alpha)n} \;\; \mbox{for}\;\; n\geq 2
\end{equation*}
a simple computation shows that 
\begin{align*}
	J_1(1)&=\sum_{n=2}^{\infty}\frac{2}{\alpha  n^2+(1-\alpha)n}+2^p\sum_{n=2}^{\infty}\frac{1}{(\alpha  n^2+(1-\alpha)n)^p}\\&\quad\quad-\sum_{n=2}^{\infty}\frac{2(-1)^{n-1}}{\alpha  n^2+(1-\alpha)n}\\&\geq  2^p\sum_{n=2}^{\infty}\frac{1}{(\alpha  n^2+(1-\alpha)n)^p}>0.
\end{align*}
Therefore, $ J_1(0)J_1(1)<0 $, and hence by the Intermediate value theorem, we conclude that $ J_1(r) $ has roots in $ (0,1) $. In order to show that $ J_1(r) $ has exactly one root in $ (0,1) $, it is sufficient to show that $ J_1 $ is strictly monotic function in $ (0,1) $. A simple computation shows that
\begin{equation*}
	J_1^{\prime}(r)=1+\sum_{n=2}^{\infty}\frac{2nr^{n-1}}{\alpha  n^2+(1-\alpha)n}+\sum_{n=2}^{\infty}\frac{np2^pr^{np-1}}{(\alpha  n^2+(1-\alpha)n)^p}>0
\end{equation*}
for all $ r\in (0,1) $, which shows that $ J_1 $ is strictly increasing function. Therefore, $ J_1(r) $ has the unique root in $ (0,1) $, say $ r_p(\alpha) $. Therefore, we have $ J_1(r_p(\alpha))=0 $ and hence from \eqref{ee-33.22}, we obtain
\begin{align}\label{e-3.2}
&r_p(\alpha)+\sum_{n=2}^{\infty}\frac{2r^n_p(\alpha)}{\alpha  n^2+(1-\alpha)n}+\sum_{n=2}^{\infty}\frac{2^pr^{np}_p(\alpha)}{(\alpha  n^2+(1-\alpha)n)^p}\\&\nonumber=1+\sum_{n=2}^{\infty}\frac{2(-1)^{n-1}}{\alpha n^2+(1-\alpha)n}.
\end{align}
To show that $ r_p(\alpha) $ is the best possible, we consider the function $ f=f_{\alpha} $ defined by \eqref{e-1.4}. It is easy to see that, the function $ f_{\alpha}\in \mathcal{W}^{0}_{\mathcal{H}}(\alpha) $ and for $ f=f_{\alpha} $, we have 
\begin{equation}\label{e-3.3}
	d(f(0), \partial f(\mathbb{D}))=1+\sum_{n=2}^{\infty}\frac{2(-1)^{n-1}}{\alpha n^2+(1-\alpha)n}.
\end{equation}
It is easy to see that 
\begin{equation*}
	r+\sum_{n=2}^{\infty}\frac{2r^n}{\alpha  n^2+(1-\alpha)n}+\sum_{n=2}^{\infty}\frac{2^pr^{np}}{(\alpha  n^2+(1-\alpha)n)^p}\leq 1+\sum_{n=2}^{\infty}\frac{2(-1)^{n-1}}{\alpha n^2+(1-\alpha)n}
\end{equation*} 
for $ r\leq r_p(\alpha) $. A simple computation using \eqref{e-3.2} and \eqref{e-3.3} for the function $ f=f_{\alpha} $ and $ r= r_p(\alpha) $,  shows that
 \begin{align*}
 	&|z|+\sum_{n=2}^{\infty}(|a_n|+|b_n|)|z|^n+\sum_{n=2}^{\infty}(|a_n|+|b_n|)^p|z|^{np}\\&=r+\sum_{n=2}^{\infty}\frac{2r^n}{\alpha  n^2+(1-\alpha)n}+\sum_{n=2}^{\infty}\frac{2^pr^{np}}{(\alpha  n^2+(1-\alpha)n)^p}\\&=r_p(\alpha)+\sum_{n=2}^{\infty}\frac{2r^n_p(\alpha)}{\alpha  n^2+(1-\alpha)n}+\sum_{n=2}^{\infty}\frac{2^pr^{np}_p(\alpha)}{(\alpha  n^2+(1-\alpha)n)^p}\\&=1+\sum_{n=2}^{\infty}\frac{2(-1)^{n-1}}{\alpha n^2+(1-\alpha)n}\\&=d(f(0),\partial f(\mathbb{D})).
 \end{align*}
Therefore, $ r_p(\alpha) $ is the best possible. This completes the proof.
\end{proof}	

\begin{proof}[\bf Proof of Corollary \ref{cor-2.4}]
		Let $ f\in \mathcal{W}^{0}_{\mathcal{H}}(\alpha)  $ be given by \eqref{e-1.1}. Then in view of Lemma \ref{lem-1.3} and Lemma \ref{lem-1.5}, and \eqref{e-3.1}, for $ |z|=r $ and $ \alpha=1/2 $, we obtain
\begin{align}\label{e-33.4}
& |z|+\sum_{n=2}^{\infty}\left(|a_n|+|b_n|\right)|z|^n+\sum_{n=2}^{\infty}\left(|a_n|+|b_n|\right)^2|z|^{2n}\\&\nonumber\leq r+\sum_{n=2}^{\infty}\frac{4r^n}{n^2+n}+\sum_{n=2}^{\infty}\left(\frac{4r^n}{n^2+n}\right)^2.
\end{align}
A simple computation shows that
\begin{align*}
		\sum_{n=2}^{\infty}\frac{4r^n}{n^2+n}&=\sum_{n=2}^{\infty}-\frac{4r^n}{n+1}+\sum_{n=2}^{\infty}\frac{4r^n}{n}\\&=-\frac{4}{r}\left(-\frac{r^2}{2}-r+\sum_{n=1}^{\infty}\frac{r^n}{n}\right)-4r+4\sum_{n=1}^{\infty}\frac{r^n}{n}\\&=-\frac{4}{r}\left(-\frac{r^2}{2}-r-\log(1-r)\right)-4r-4\log(1-r)\\&=4-2r+\frac{4}{r}(1-r)\log (1-r).
\end{align*}
On the other hand, it is easy to see that
\begin{align*}
	\sum_{n=2}^{\infty}\left(\frac{4r^n}{n^2+n}\right)^2&=\sum_{n=2}^{\infty}-\frac{32r^{2n}}{n}+\sum_{n=2}^{\infty}\frac{16r^{2n}}{n^2}+\sum_{n=2}^{\infty}\frac{16r^{2n}}{(n+1)^2}+\sum_{n=2}^{\infty}\frac{32r^{2n}}{n+1}.
\end{align*}
Using the defintion of dilogarithm function $ {\rm Li}_2 $, a simple computation shows that
	\begin{equation}\label{e-33.5}
		\begin{cases}
		\displaystyle\sum_{n=2}^{\infty}-\frac{32r^{2n}}{n}=32r^2+32\log(1-r^2),\vspace{1.8mm}\\
		\displaystyle\sum_{n=2}^{\infty}\frac{16r^{2n}}{n^2}=-16r^2+16{\rm Li}_2(r^2),\vspace{1.8mm}\\
		\displaystyle\sum_{n=2}^{\infty}\frac{16r^{2n}}{(n+1)^2}=\frac{16}{r^2}\left(-\frac{r^4}{4}-r^2+{\rm Li}_2\left(r^2\right)\right),\vspace{1.8mm}\\
		\displaystyle\sum_{n=2}^{\infty}\frac{32r^{2n}}{n+1}=\frac{32}{r^2}\left(-\frac{r^4}{2}-r^2+\left(-\log\left(1-r^2\right)\right)\right).
	\end{cases}
	\end{equation}
Therefore, using \eqref{e-33.5}, we obtain
\begin{align*}
	\sum_{n=2}^{\infty}\left(\frac{4r^n}{n^2+n}\right)^2&=\left(32r^2+32\log(1-r^2)\right)+\left(-16r^2+16{\rm Li}_2(r^2)\right)\\&\quad+\frac{16}{r^2}\left(-\frac{r^4}{4}-r^2+{\rm Li}_2\left(r^2\right)\right)+\frac{32}{r^2}\left(-\frac{r^4}{2}-r^2+\left(-\log\left(1-r^2\right)\right)\right)\\&=-4r^2+32\log (1-r^2)+16{\rm Li}_2(r^2)-48-\frac{32\log(1-r^2)}{r^2}\\&\quad\quad+\frac{16{\rm Li}_2\left(r^2\right)}{r^2}\\&=-4r^2-48+\frac{16\left(1+r^2\right){\rm Li}_2\left(r^2\right)}{r^2}-\frac{32}{r^2}\left(1-r^2\right)\log\left(1-r^2\right).
\end{align*}
Hence from \eqref{e-3.1} and \eqref{e-33.4}, we obtain
\begin{align*}
	&|z|+\sum_{n=2}^{\infty}\left(|a_n|+|b_n|\right)|z|^n+\sum_{n=2}^{\infty}\left(|a_n|+|b_n|\right)^2|z|^{2n}\\&=-4r^2-r-44+\frac{16}{r^2}\left(\left(1+r^2\right){\rm Li}_2\left(r^2\right)\right)+\frac{4}{r}\left((1-r)\log(1-r)\right)\\&\quad\quad-\frac{32}{r^2}\left(\left(1-r^2\right)\log\left(1-r^2\right)\right).
\end{align*} 
A simple computation shows that
\begin{align*}
	&-4r^2-r-44+\frac{16}{r^2}\left(\left(1+r^2\right){\rm Li}_2\left(r^2\right)\right)+\frac{4}{r}\left((1-r)\log(1-r)\right)\\&\quad\quad-\frac{32}{r^2}\left(\left(1-r^2\right)\log\left(1-r^2\right)\right)\leq 1+2(-3+4\log 2)
\end{align*} 
for $ r\leq r_2(0.5) $, where $ r_2(0.5)\approx 0.399085 $ is a root of $ G(r)=0 $ in $ (0,1) $ where $ G : [0,1]\rightarrow\mathbb{R} $ be defined by
\begin{align*}
	G(r):&=-4r^2-r+\frac{16}{r^2}\left(\left(1+r^2\right){\rm Li}_2\left(r^2\right)\right)+\frac{4}{r}\left((1-r)\log(1-r)\right)\\&\quad\quad-\frac{32}{r^2}\left(\left(1-r^2\right)\log\left(1-r^2\right)\right)-39-8\log 2.
\end{align*}
By following the argument as in the proof of Theorem \ref{th-2.1}, we can show that the function $ G(r) $ has the unqiue root $ r_2(0.5)\approx 0.399085 $ in $ (0,1) $. Thus, we have
\begin{align*}
	&-4r^2_2(0.5)-r_2(0.5)+\frac{16\left(1+r^2_2(0.5)\right){\rm Li}_2\left(r^2_2(0.5)\right)}{r^2_2(0.5)}+\frac{4(1-r_2(0.5))\log(1-r_2(0.5))}{r}\\&\quad\quad-\frac{32\left(1-r^2_2(0.5)\right)\log\left(1-r^2_2(0.5)\right)}{r^2_2(0.5)}-39-8\log 2=0
\end{align*}
which is equivalent to
\begin{align}\label{e-33.6}
	&-4r^2_2(0.5)+\frac{16\left(1+r^2_2(0.5)\right){\rm Li}_2\left(r^2_2(0.5)\right)}{r^2_2(0.5)}+\frac{4(1-r_2(0.5))\log(1-r_2(0.5))}{r}\\&\nonumber\quad\quad-r_2(0.5)-44-\frac{32\left(1-r^2_2(0.5)\right)\log\left(1-r^2_2(0.5)\right)}{r^2_2(0.5)}= 1+2(-3+4\log 2).
\end{align}
To show that $ r_2(0.5) $ is the best possible, we consider the function $ f=f_{1/2} $ defined by
	\begin{equation}\label{e-33.7}
		f_{1/2}(z)=z+\sum_{n=2}^{\infty}\frac{4z^n}{n^2+n}.
	\end{equation}
 In view of \eqref{e-3.1}, it is easy to see that
\begin{equation}\label{e-33.8}
	d(f(0),\partial f(\mathbb{D}))=1+\sum_{n=2}^{\infty}\frac{4(-1)^{n-1}}{n^2+n}=1+2 (-3+4\log 2).
\end{equation}
A simple computation using \eqref{e-33.4}, \eqref{e-33.6} and \eqref{e-33.8}, for $ f=f_{1/2} $ and $ r=r_2(0.5) $ shows that
\begin{align*}
	&|z|+\sum_{n=2}^{\infty}\left(|a_n|+|b_n|\right)|z|^n+\sum_{n=2}^{\infty}\left(|a_n|+|b_n|\right)^2|z|^{2n}\\&=-4r^2_2(0.5)+\frac{16\left(1+r^2_2(0.5)\right){\rm Li}_2\left(r^2_2(0.5)\right)}{r^2_2(0.5)}+\frac{4(1-r_2(0.5))\log(1-r_2(0.5))}{r}\\&\nonumber\quad\quad-r_2(0.5)-44-\frac{32\left(1-r^2_2(0.5)\right)\log\left(1-r^2_2(0.5)\right)}{r^2_2(0.5)}\\&= 1+2(-3+4\log 2)\\&=d(f_{1/2}(0),\partial f_{1/2}(\mathbb{D})).
\end{align*}
Therefore, $ r_2(0.5) $ is the best possible. This completes the proof.
\end{proof}

\begin{proof}[\bf Proof of Corollary \ref{cor-2.5}]
	Let $ f\in \mathcal{W}^{0}_{\mathcal{H}}(\alpha)  $ be given by \eqref{e-1.1}. Then in view of Lemma \ref{lem-1.3} and Lemma \ref{lem-1.5}, and \eqref{e-3.1}, for $ |z|=r $ and $ \alpha=1/2 $, we obtain
	\begin{align}\label{e-33.9}
		&|z|+\sum_{n=2}^{\infty}\left(|a_n|+|b_n|\right)^2|z|^{2n}\\&\nonumber\leq r+\sum_{n=2}^{\infty}\frac{16r^{2n}}{(n^2+n)^2}\\&\nonumber=r-4r^2-48+\frac{16\left(1+r^2\right){\rm Li}_2\left(r^2\right)}{r^2}-\frac{32}{r^2}\left(1-r^2\right)\log\left(1-r^2\right).
	\end{align}
A simple computation shows that
\begin{align*}
		&r-4r^2-48+\frac{16\left(1+r^2\right){\rm Li}_2\left(r^2\right)}{r^2}-\frac{32}{r^2}\left(1-r^2\right)\log\left(1-r^2\right)\\&\leq 1+2(-3+4\log 2)
\end{align*}
for $ r\leq r^*(0.5) $, where $ r^*(0.5)\approx 0.512331 $ is a root of $ H(r)=0 $ in $ (0,1) $ where $ H : [0,1]\rightarrow\mathbb{R} $ be defined by
\begin{align*}
	H(r) : = r-4r^2-43+\frac{16\left(1+r^2\right){\rm Li}_2\left(r^2\right)}{r^2}-\frac{32}{r^2}\left(1-r^2\right)\log\left(1-r^2\right)-8\log 2.
\end{align*}
By using the same urgument as in the proof of Theorem \ref{th-2.1}, we can show that $ H(r) $ has the unique root $ r^*(0.5)\approx 0.512331 $ in $ (0,1) $. Therefore, we have $ H(r^*(0.5))=0$. That is,
\begin{align*}
&r^*(0.5)-4\left(r^*(0.5)\right)^2-43+\frac{16\left(1+\left(r^*(0.5)\right)^2\right){\rm Li}_2\left(\left(r^*(0.5)\right)^2\right)}{\left(r^*(0.5)\right)^2}\\&\quad\quad-\frac{32}{\left(r^*(0.5)\right)^2}\left(1-\left(r^*(0.5)\right)^2\right)\log\left(1-\left(r^*(0.5)\right)^2\right)-8\log 2=0
\end{align*}
which is equivalent to
\begin{align}\label{e-33.10}
	&r^*(0.5)-4\left(r^*(0.5)\right)^2-48+\frac{16\left(1+\left(r^*(0.5)\right)^2\right){\rm Li}_2\left(\left(r^*(0.5)\right)^2\right)}{\left(r^*(0.5)\right)^2}\\&\quad\quad\nonumber-\frac{32}{\left(r^*(0.5)\right)^2}\left(1-\left(r^*(0.5)\right)^2\right)\log\left(1-\left(r^*(0.5)\right)^2\right)\\&\nonumber= 1+2(-3+4\log 2).
\end{align}
To show that $ r^*(0.5) $ is the best possible, we consider the function $ f=f_{1/2} $ defined by \eqref{e-33.7}. Then a simple computation using \eqref{e-33.8}, \eqref{e-33.9} and \eqref{e-33.10} for $ f=f_{1/2} $ and $ r=r^*(0.5) $ shows that
\begin{align*}
&|z|+\sum_{n=2}^{\infty}\left(|a_n|+|b_n|\right)^2|z|^{2n}\\&= r^*(0.5)+\sum_{n=2}^{\infty}\frac{16(r^*(0.5))^{2n}}{(n^2+n)^2}\\&=r^*(0.5)-4\left(r^*(0.5)\right)^2-48+\frac{16\left(1+\left(r^*(0.5)\right)^2\right){\rm Li}_2\left(\left(r^*(0.5)\right)^2\right)}{\left(r^*(0.5)\right)^2}\\&\quad\quad\nonumber-\frac{32}{\left(r^*(0.5)\right)^2}\left(1-\left(r^*(0.5)\right)^2\right)\log\left(1-\left(r^*(0.5)\right)^2\right)\\&= 1+2(-3+4\log 2)\\&=d(f_{1/2}(0), \partial f_{1/2}(\mathbb{D})).
\end{align*}
Therefore, $ r^*(0.5)\approx 0.512331 $ is the best possible. This completes the proof.
\end{proof}
\begin{proof}[\bf Proof of Theorem \ref{th-2.6}]
		Let $ f\in \mathcal{W}^{0}_{\mathcal{H}}(\alpha)  $ be given by \eqref{e-1.1}. Then in view of Lemma \ref{lem-1.3} and Lemma \ref{lem-1.5}, and \eqref{e-3.1}, for $ |z|=r $, we obtain
		\begin{align}\label{e-3.4}
			&|f(z)|+\sum_{n=2}^{\infty}\left(|a_n|+|b_n|\right)|z|^n\\&\leq \nonumber |z|+\sum_{n=2}^{\infty}\frac{2|z|^n}{\alpha n^2+(1-\alpha)n}+\sum_{n=2}^{\infty}\frac{2|z|^n}{\alpha n^2+(1-\alpha)n}\\&=\nonumber r+\sum_{n=2}^{\infty}\frac{4r^n}{\alpha n^2+(1-\alpha)n}.
		\end{align}
	A simple computation shows that 
	\begin{equation*}
		r+\sum_{n=2}^{\infty}\frac{4r^n}{\alpha n^2+(1-\alpha)n}\leq 1+\sum_{n=2}^{\infty}\frac{2(-1)^{n-1}}{\alpha n^2+(1-\alpha)n}.
	\end{equation*}
for $ r\leq r_*(\alpha) $, where $ r_*(\alpha) $ is a root of $ J_2(r)=0 $ in $ (0,1) $, where $ J_2 : [0,1]\rightarrow\mathbb{R} $ be defined by 
 \begin{equation*}
 	J_2(r):=	r+\sum_{n=2}^{\infty}\frac{4r^n}{\alpha n^2+(1-\alpha)n}-1-\sum_{n=2}^{\infty}\frac{2(-1)^{n-1}}{\alpha n^2+(1-\alpha)n}.
 \end{equation*}
It is easy to see that $ J_2(r) $ is conitnuous on $ [0,1] $ and differentiable on $ (0,1) $. Since 
\begin{equation*}
	J_2(0)=-1-\sum_{n=2}^{\infty}\frac{2(-1)^{n-1}}{\alpha n^2+(1-\alpha)n}<0
\end{equation*}
and 
\begin{align*}
	J_2(1)&=\sum_{n=2}^{\infty}\frac{4}{\alpha n^2+(1-\alpha)n}-\sum_{n=2}^{\infty}\frac{2(-1)^{n-1}}{\alpha n^2+(1-\alpha)n}\\&=\sum_{n=2}^{\infty}\frac{4-2(-1)^{n-1}}{\alpha n^2+(1-\alpha)n}>0,
\end{align*}
by the Intermediate value theorem, $ J_2 $ has a root in $ (0,1) $. Since 
\begin{equation}\label{e-3.5}
	J_2^{\prime}(r)=1+\sum_{n=2}^{\infty}\frac{4nr^{n-1}}{\alpha n^2+(1-\alpha)n}>0\;\; \mbox{for}\;\; n\geq 2,
\end{equation}
the function $ J_2(r) $ is strictly increasing function in $ (0,1) $. Therefore, the function $ J_2 $ has the unique root, say $ r_*(\alpha) $, in $ (0,1) $. Therefore, $ J_2(r_*(\alpha))=0 $ which is equivalent to
\begin{equation}\label{e-3.6}
		r_*(\alpha)+\sum_{n=2}^{\infty}\frac{4r_*^n(\alpha)}{\alpha n^2+(1-\alpha)n}=1+\sum_{n=2}^{\infty}\frac{2(-1)^{n-1}}{\alpha n^2+(1-\alpha)n}.
\end{equation}
To show that $ r_*(\alpha) $ is the best possible, we consider the function $ f=f_{\alpha} $ defined by \eqref{e-1.4}. Clearly, the function $ f_{\alpha} $ belongs to the class $ \mathcal{W}^0_{\mathcal{H}}(\alpha) $. A simple computation using \eqref{e-3.3}, \eqref{e-3.5} and \eqref{e-3.6}, for $ f=f_{\alpha} $ and $ r=r_*(\alpha) $ shows that 
\begin{align*}
	|f_{\alpha}(z)|+\sum_{n=2}^{\infty}\left(|a_n|+|b_n|\right)|z|^n&= r_*(\alpha)+\sum_{n=2}^{\infty}\frac{4r_*^n(\alpha)}{\alpha n^2+(1-\alpha)n}\\&=1+\sum_{n=2}^{\infty}\frac{2(-1)^{n-1}}{\alpha n^2+(1-\alpha)n}\\&=d(f(0), \partial f(\mathbb{D})).
\end{align*}
This shows that $ r_*(\alpha) $ is the best possible.
\end{proof}	
\begin{proof}[\bf Proof of Corollary \ref{cor-2.10}]
	Let $ f\in \mathcal{W}^{0}_{\mathcal{H}}(\alpha)  $ be given by \eqref{e-1.1}. Then for $ \alpha=1/2 $, using Lemma \ref{lem-1.3} and Lemma \ref{lem-1.5} and \eqref{e-3.4}, we obtain 
\begin{equation}\label{e-3.7}
	|f(z)|+\sum_{n=2}^{\infty}\left(|a_n|+|b_n|\right)|z|^n\leq r+4\sum_{n=2}^{\infty}\frac{2r^n}{n^2+n}.	
\end{equation}
A simple computation shows that
\begin{align*}
	\sum_{n=2}^{\infty}\frac{2r^n}{n^2+n}&=\sum_{n=2}^{\infty}-\frac{2r^n}{n+1}+\sum_{n=2}^{\infty}\frac{2r^n}{n}\\&=\sum_{n=3}^{\infty}-\frac{2r^n}{nr}+\sum_{n=2}^{\infty}\frac{2r^n}{n}\\&=-\frac{2}{r}\left(-r-\frac{r^2}{2}+\sum_{n=1}^{\infty}\frac{r^n}{n}\right)-2r+2\sum_{n=1}^{\infty}\frac{r^n}{n}\\&=-\frac{2}{r}\left(-r-\frac{r^2}{2}-\log (1-r)\right)-2r-2\log (1-r)\\&=-r-2\log (1-r)+2+\frac{2}{r}\log (1-r).
\end{align*}
On the other hand, a simple computation shows that
\begin{equation*}
	\sum_{n=2}^{\infty}\frac{2(-1)^{n-1}}{n^2+n}=-3+4\log 2.
\end{equation*}
In view of \eqref{e-3.1} and \eqref{e-3.7}, we obtain 
\begin{equation*}
		|f(z)|+\sum_{n=2}^{\infty}\left(|a_n|+|b_n|\right)|z|^n\leq 	r+4\left(-r-2\log (1-r)+2+\frac{2}{r}\log (1-r)\right).
\end{equation*}
A simple computation shows that
\begin{equation*}
	r+4\left(-r-2\log (1-r)+2+\frac{2}{r}\log (1-r)\right)\leq 1+2 (-3+4\log 2)
\end{equation*}
for $ r\leq  r_*(0.5)$, where $ r_*(0.5) $ is a root of $ H_1(r)=0 $ $ (0,1) $ where $ H_1 : [0,1]\rightarrow\mathbb{R} $ be defined by
\begin{align*}
		H_1(r):&=r+4\left(-r-2\log (1-r)+2+\frac{2}{r}\log (1-r)\right)- 1-2 (-3+4\log 2)\\&=\frac{8}{r}(1-r)\log(1-r)-3r+13-8\log 2.
\end{align*}
By following the same argument as in the proof of Theorem \ref{th-2.6}, we can show that $ H_1(r) $ has the unique root $ r_*(0.5)\approx 0.347966 $. This shows that $ H_1(r_*(0.5))=0 $. That is,
\begin{equation*}
	\frac{8}{r_*(0.5)}(1-r_*(0.5))\log(1-r_*(0.5))-3r_*(0.5)+13-8\log 2
\end{equation*}
which is equivalent to 
\begin{align}
	\label{e-3.8} &
	r_*(0.5)+4\left(-r_*(0.5)-2\log (1-r_*(0.5))+2+\frac{2}{r_*(0.5)}\log (1-r_*(0.5))\right)\\&=\nonumber 1+2 (-3+4\log 2).
\end{align}
To prove that $ r_*(0.5) $ is the best possible, we consider the function $ f=f_{1/2} $ given by \eqref{e-33.7}. Further, a simple computation using \eqref{e-3.7}, \eqref{e-33.8} and \eqref{e-3.8} for $ f=f_{1/2} $ and $ r=r_*(0.5) $  shows that
\begin{align*}
	&|f_{1/2}(z)|+\sum_{n=2}^{\infty}\left(|a_n|+|b_n|\right)|z|^n\\&= r_*(0.5)+4\sum_{n=2}^{\infty}\frac{2r^n_*(0.5)}{n^2+n}\\&=	r_*(0.5)+4\left(-r_*(0.5)-2\log (1-r_*(0.5))+2+\frac{2}{r_*(0.5)}\log (1-r_*(0.5))\right)\\&=1+2 (-3+4\log 2)\\&=d(f_{1/2}(0),\partial f_{1/2}(\mathbb{D})).
\end{align*}
Therefore, $ r_*(0.5)\approx 0.347966 $ is the best possible. This completes the proof.
\end{proof}	
\begin{proof}[\bf Proof of Theorem \ref{th-2.11}]
Let $ f\in \mathcal{W}^{0}_{\mathcal{H}}(\alpha)  $ be given by \eqref{e-1.1}. Then in view of Lemma \ref{lem-1.3} and Lemma \ref{lem-1.5}, and \eqref{e-3.1}, for $ |z|=r $, we obtain 
\begin{align}\label{e-3.16}
	|f(z^m)|+\sum_{n=N}^{\infty}\left(|a_n|+|b_n|\right)|z|^n&\leq |z|^m+\sum_{n=2}^{\infty}\frac{2|z|^{mn}}{\alpha n^2+(1-\alpha)n}+\sum_{n=N}^{\infty}\frac{2|z|^n}{\alpha n^2+(1-\alpha)n}\\&\nonumber=r^m+\sum_{n=2}^{\infty}\frac{2r^{mn}}{\alpha n^2+(1-\alpha)n}+\sum_{n=N}^{\infty}\frac{2r^n}{\alpha n^2+(1-\alpha)n}.
\end{align}
A simple computation shows that
	\begin{align*}
		r^m+\sum_{n=2}^{\infty}\frac{2r^{mn}}{\alpha n^2+(1-\alpha)n}+\sum_{n=N}^{\infty}\frac{2r^n}{\alpha n^2+(1-\alpha)n}\leq 1+\sum_{n=2}^{\infty}\frac{2(-1)^{n-1}}{\alpha n^2+(1-\alpha)n}
	\end{align*}
for $ r\leq  R_{m,N}(\alpha)$, where $ R_{m,N}(\alpha) $ is a root of $ J_3(r)=0 $ in $ (0,1) $, where $ J_3 : [0,1]\rightarrow\mathbb{R} $ be defined by 
\begin{align*}
	J_3(r):&=r^m+\sum_{n=2}^{\infty}\frac{2r^{mn}}{\alpha n^2+(1-\alpha)n}+\sum_{n=N}^{\infty}\frac{2r^n}{\alpha n^2+(1-\alpha)n}\\&\quad\quad -1-\sum_{n=2}^{\infty}\frac{2(-1)^{n-1}}{\alpha n^2+(1-\alpha)n}.
\end{align*}
A simple computation shows that $ J_3(0)J_3(1)<0 $. Further, $ J_3^{\prime}(r) $ is strictly increasing evident from the fact that
\begin{align*}
	J_3^{\prime}(r)=mr^{m-1}+\sum_{n=2}^{\infty}\frac{2mnr^{mn-1}}{\alpha n^2+(1-\alpha)n}+\sum_{n=N}^{\infty}\frac{2nr^{n-1}}{\alpha n^2+(1-\alpha)n}>0\;\; \mbox{for all}\;\; r\in (0,1).
\end{align*}
 The function $ J_3 $ being differentiable and stirctly increasing on $ (0,1) $, Intermadiate value theorem asserts that $ J_3(r) $ has the unique root in $ (0,1) $, say $ R_{m,N}(\alpha) $. Therefore, we have  $ J_3(R_{m,N}(\alpha))=0 $ which is equivalent to
\begin{align}\label{e-3.17}
	&R^m_{m,N}(\alpha)+\sum_{n=2}^{\infty}\frac{2R^{mn}_{m,N}(\alpha)}{\alpha n^2+(1-\alpha)n}+\sum_{n=N}^{\infty}\frac{2R^n_{m,N}(\alpha)}{\alpha n^2+(1-\alpha)n}\\&\nonumber= 1+\sum_{n=2}^{\infty}\frac{2(-1)^{n-1}}{\alpha n^2+(1-\alpha)n}.
\end{align}
To show that $ R_{m,N}(\alpha) $ is the best possible, we consider the function $ f=f_{\alpha} $ defined by \eqref{e-1.4}. In view of \eqref{e-3.3}, \eqref{e-3.17} and \eqref{e-3.16}, for $ f=f_{\alpha} $ and $ z= R_{m,N}(\alpha)$, we obtain
\begin{align*}
&|f_{\alpha}(z^m)|+\sum_{n=N}^{\infty}\left(|a_n|+|b_n|\right)|z|^n\\&=	R^m_{m,N}(\alpha)+\sum_{n=2}^{\infty}\frac{2R^{mn}_{m,N}(\alpha)}{\alpha n^2+(1-\alpha)n}+\sum_{n=N}^{\infty}\frac{2R^n_{m,N}(\alpha)}{\alpha n^2+(1-\alpha)n}\\&= 1+\sum_{n=2}^{\infty}\frac{2(-1)^{n-1}}{\alpha n^2+(1-\alpha)n}\\&= d(f_{\alpha}(0),\partial f_{\alpha}(\mathbb{D})).
\end{align*}
This shows that $ R_{m,N}(\alpha) $ is the best possible. This completes the proof.
\end{proof}	
\begin{proof}[\bf Proof of Theorem \ref{th-2.13}]
	For $ f\in\mathcal{W}^0_{\mathcal{H}}(\alpha) $, the Jacobian of $ f $ is denoted by $ J_f $ and is defined by
	\begin{equation*}
		J_f(z)=|f_{z}(z)|^2-|f_{\bar{z}}(z)|^2=|h^{\prime}(z)|^2-|g^{\prime}(z)|^2\;\; \mbox{for}\;\; z\in\mathbb{D}.
	\end{equation*}
It is weel-known that (see \cite[p.113]{Duren-2004}) the area of the disk $ \mathbb{D}_r:=\{z\in\mathbb{C} : |z|<r\} $ under the harmonic map $ f=h+\bar{g} $ is 
\begin{align}\label{e-3.18}
	S_r=\iint\limits_{\mathbb{D}_r}J_f(z)dxdy=\iint\limits_{\mathbb{D}_r}\left(|{h^{\prime}(z)}|^2-|{g^{\prime}(z)}|^2\right)dxdy.
\end{align}
Therefore, by using polar coordinates, we obtain 
\begin{align}\label{e-3.19}
	\iint\limits_{\mathbb{D}_r}|h^{\prime}(z)|^2dxdy&=\int_{0}^{r}\int_{0}^{2\pi}|h^{\prime}(\rho e^{i\theta})|^2\rho d\theta d\rho\\&=\nonumber\int_{0}^{r}\int_{0}^{2\pi}\rho\left(\sum_{n=1}^{\infty}n\;a_n\rho^{n-1}e^{i(n-1)\theta}\right)\left(\sum_{n=1}^{\infty}n\;\bar{a}_n\rho^{n-1}e^{-i(n-1)\theta}\right)d\theta d\rho\\&=\nonumber \int_{0}^{r}\left(\sum_{n=1}^{\infty}2\pi n^2 |a_n|^2\rho^{2n-1}\right)d\rho\\&\nonumber=\sum_{n=1}^{\infty}2\pi n^2|a_n|^2\frac{r^{2n}}{2n}\\&\nonumber=\pi \sum_{n=1}^{\infty}n|a_n|^2r^{2n}.
\end{align}
Similarly, for $ g(z)=\sum_{n=2}^{\infty}b_nz^n $, a simple computation shows that
\begin{align}\label{e-3.20}
	\iint\limits_{\mathbb{D}_r}&|g^{\prime}(z)|^2dxdy=\pi\sum_{n=2}^{\infty}n|b_n|^2r^{2n}.
\end{align}
In view of Lemma \ref{lem-1.3}, \eqref{e-3.18}, \eqref{e-3.19} and  \eqref{e-3.20}, we obtain
\begin{align}\label{e-3.21}
	\frac{S_r}{\pi}&=\frac{1}{\pi}\iint\limits_{\mathbb{D}_r}\left(|h^{\prime}(z)|^2-|g^{\prime}(z)|^2\right)dxdy\\&\nonumber=r^2+\sum_{n=2}^{\infty}n\left(|a_n|^2-|b_n|^2\right)r^{2n}\\&\nonumber=r^2+\sum_{n=2}^{\infty}n(|a_n|+|b_n|)(|a_n|-|b_n|)r^{2n}\\&\nonumber=r^2+\sum_{n=2}^{\infty}\frac{4nr^{2n}}{\left(\alpha n^2+(1-\alpha)n\right)^2}.
\end{align}
(i) In view of Lemma \ref{lem-1.3} and Lemma \ref{lem-1.5} and  \eqref{e-3.21} for $ |z|=r $, we obtain
\begin{align}\label{e-3.22}
	&|z|+\sum_{n=2}^{\infty}\left(|a_n|+|b_n|\right)|z|^n+\frac{S_r}{\pi}\\&\nonumber\leq r+\sum_{n=2}^{\infty}\frac{2r^n}{\alpha n^2+(1-\alpha)n}+r^2+\sum_{n=2}^{\infty}\frac{4nr^{2n}}{\left(\alpha n^2+(1-\alpha)n\right)^2}.
\end{align}
It is easy to see that
\begin{align*}
	 r^2+r+\sum_{n=2}^{\infty}\frac{2r^n}{\alpha n^2+(1-\alpha)n}+\sum_{n=2}^{\infty}\frac{4nr^{2n}}{\left(\alpha n^2+(1-\alpha)n\right)^2}\leq 1+\sum_{n=2}^{\infty}\frac{2(-1)^{n-1}}{\alpha n^2+(1-\alpha)n}.
\end{align*}
for $ r\leq  r_f(\alpha)$, where $ r_f(\alpha) $ is the root of $ J_4(r)=0 $, where $ J_4 : [0,1]\rightarrow\mathbb{R} $ is defined by 
\begin{align*}
	J_4(r):&= r^2+r+\sum_{n=2}^{\infty}\frac{2r^n}{\alpha n^2+(1-\alpha)n}+\sum_{n=2}^{\infty}\frac{4nr^{2n}}{\left(\alpha n^2+(1-\alpha)n\right)^2}\\&\quad\quad- 1-\sum_{n=2}^{\infty}\frac{2(-1)^{n-1}}{\alpha n^2+(1-\alpha)n}.
\end{align*}
It is not difficult to show that $ J_4(0)J_4(1)<0 $ and $ J^{\prime}_4(r)>0 $ for $ r\in (0,1) $. Then by the Intermediate value theorem, the function $ J_4 $ has the unique root in $ (0,1) $, say $ r_f(\alpha) $. Therefore, we have
\begin{align}\label{e-3.23}
	&r^2_f(\alpha)+r_f(\alpha)+\sum_{n=2}^{\infty}\frac{2r^n_f(\alpha)}{\alpha n^2+(1-\alpha)n}+\sum_{n=2}^{\infty}\frac{4nr^{2n}_f(\alpha)}{\left(\alpha n^2+(1-\alpha)n\right)^2}\\&\nonumber= 1+\sum_{n=2}^{\infty}\frac{2(-1)^{n-1}}{\alpha n^2+(1-\alpha)n}.
\end{align}
To show that $ r_f(\alpha) $ is the best possible, we consider the function $ f=f_{\alpha} $ given by \eqref{e-1.4}. A simple computation using \eqref{e-3.3}, \eqref{e-3.22} and \eqref{e-3.23} for $ f=f_{\alpha} $ and $ r=r_f(\alpha) $ shows that 
\begin{align*}
	&|z|+\sum_{n=2}^{\infty}\left(|a_n|+|b_n|\right)|z|^n+\frac{S_{r_f(\alpha)}}{\pi}\\&=r^2_f(\alpha)+r_f(\alpha)+\sum_{n=2}^{\infty}\frac{2r^n_f(\alpha)}{\alpha n^2+(1-\alpha)n}+\sum_{n=2}^{\infty}\frac{4nr^{2n}_f(\alpha)}{\left(\alpha n^2+(1-\alpha)n\right)^2}\\&= 1+\sum_{n=2}^{\infty}\frac{2(-1)^{n-1}}{\alpha n^2+(1-\alpha)n}\\&=d(f_{\alpha}(0), \partial f_{\alpha}(\mathbb{D})). 
\end{align*}
Therefore, $ r_f(\alpha) $ is the best possible. This complete the proof of (i).\vspace{1.5mm}

\noindent (ii) In view of Lemma \ref{lem-1.3} and Lemma \ref{lem-1.5} and \eqref{e-3.21} for $ |z|=r $, we obtain
\begin{align}\label{e-3.24}
	&|f(z)|^2+\sum_{n=2}^{\infty}\left(|a_n|+|b_n|\right)|z|^n+\left(\frac{S_r}{\pi}\right)^2\\&\leq \nonumber \left(|z|+\sum_{n=2}^{\infty}(|a_n|+|b_n|)|z|^n\right)^2+\sum_{n=2}^{\infty}(|a_n|+|b_n|)|z|^n\\&\nonumber\quad\quad+\left(r^2+\sum_{n=2}^{\infty}\frac{4nr^{2n}}{\left(\alpha n^2+(1-\alpha)n\right)^2}\right)^2\\&\leq \nonumber \left(r+\sum_{n=2}^{\infty}\frac{2r^n}{\alpha n^2+(1-\alpha)n}\right)^2+\sum_{n=2}^{\infty}\frac{2r^n}{\alpha n^2+(1-\alpha)n}\\&\nonumber\quad\quad+\left(r^2+\sum_{n=2}^{\infty}\frac{4nr^{2n}}{\left(\alpha n^2+(1-\alpha)n\right)^2}\right)^2.
\end{align}
A simple computation shows that
\begin{align*}
	&\left(r+\sum_{n=2}^{\infty}\frac{2r^n}{\alpha n^2+(1-\alpha)n}\right)^2+\sum_{n=2}^{\infty}\frac{2r^n}{\alpha n^2+(1-\alpha)n}\\&\nonumber\quad\quad+\left(r^2+\sum_{n=2}^{\infty}\frac{4nr^{2n}}{\left(\alpha n^2+(1-\alpha)n\right)^2}\right)^2\leq 1+\sum_{n=2}^{\infty}\frac{2(-1)^{n-1}}{\alpha n^2+(1-\alpha)n}.
\end{align*}
for $ r\leq  r^*_f(\alpha)$, where $ r^*_f(\alpha) $ is a root of $ J_5(r)=0 $ in $ (0,1) $, where $ J_5 : [0,1]\rightarrow\mathbb{R} $ is defined by
\begin{align*}
	J_5(r)&=\left(r+\sum_{n=2}^{\infty}\frac{2r^n}{\alpha n^2+(1-\alpha)n}\right)^2+\sum_{n=2}^{\infty}\frac{2r^n}{\alpha n^2+(1-\alpha)n}\\&\nonumber\quad\quad+\left(r^2+\sum_{n=2}^{\infty}\frac{4nr^{2n}}{\left(\alpha n^2+(1-\alpha)n\right)^2}\right)^2- 1-\sum_{n=2}^{\infty}\frac{2(-1)^{n-1}}{\alpha n^2+(1-\alpha)n}.
\end{align*}
Using the same argument as in the proof of the Theorem \ref{th-2.1}, we can easily show that $ J_5(r) $ has the unique root, say $ r^*_f(\alpha) $. Therefore, we have
\begin{align}\label{e-3.25}
	&\left(r^*_f(\alpha)+\sum_{n=2}^{\infty}\frac{2\left(r^*_f(\alpha)\right)^n}{\alpha n^2+(1-\alpha)n}\right)^2+\sum_{n=2}^{\infty}\frac{2\left(r^*_f(\alpha)\right)^n}{\alpha n^2+(1-\alpha)n}\\&\nonumber\quad\quad+\left(\left(r^*_f(\alpha)\right)^2+\sum_{n=2}^{\infty}\frac{4n\left(r^*_f(\alpha)\right)^{2n}}{\left(\alpha n^2+(1-\alpha)n\right)^2}\right)^2= 1+\sum_{n=2}^{\infty}\frac{2(-1)^{n-1}}{\alpha n^2+(1-\alpha)n}.
\end{align}
To show $ r^*_f(\alpha) $ is the best possible, we consider the function $ f=f_{\alpha} $ given in \eqref{e-1.4}. In view of \eqref{e-3.3}, \eqref{e-3.24} and \eqref{e-3.25}, for $ f=f_{\alpha} $ and $ r=r^*_f(\alpha) $, we obtain
\begin{align*}
		&|f_{\alpha}(z)|^2+\sum_{n=2}^{\infty}\left(|a_n|+|b_n|\right)|z|^n+\left(\frac{S_{r^*_f(\alpha)}}{\pi}\right)^2\\&=	\left(r^*_f(\alpha)+\sum_{n=2}^{\infty}\frac{2\left(r^*_f(\alpha)\right)^n}{\alpha n^2+(1-\alpha)n}\right)^2+\sum_{n=2}^{\infty}\frac{2\left(r^*_f(\alpha)\right)^n}{\alpha n^2+(1-\alpha)n}\\&\nonumber\quad\quad+\left(\left(r^*_f(\alpha)\right)^2+\sum_{n=2}^{\infty}\frac{4n\left(r^*_f(\alpha)\right)^{2n}}{\left(\alpha n^2+(1-\alpha)n\right)^2}\right)^2\\&= 1+\sum_{n=2}^{\infty}\frac{2(-1)^{n-1}}{\alpha n^2+(1-\alpha)n}\\&=d(f_{\alpha}(0),\partial f_{\alpha}(\mathbb{D})).
\end{align*} 
This shows that $ r^*_f(\alpha) $ is the best possible. This completes the proof of (ii).
\end{proof}	
\begin{proof}[\bf Proof of Corollary \ref{cor-2.14}]
	(i) Let $f\in \mathcal{W}^0_{\mathcal{H}}(\alpha) $, then for $ |z|=r $ and $ \alpha=1/2 $, using Lemma \ref{lem-1.3}, we obtain
	\begin{align}\label{e-33.26}
		|z|+\sum_{n=2}^{\infty}\left(|a_n|+|b_n|\right)|z|^n+\frac{S_r}{\pi}\leq 	r^2+r+\sum_{n=2}^{\infty}\frac{4r^n}{n^2+n}+\sum_{n=2}^{\infty}\frac{16nr^{2n}}{(n^2+n)^2}.
	\end{align}
 A simple computation shows that
 \begin{equation*}
 \sum_{n=2}^{\infty}\frac{4r^n}{n^2+n}=4-2r+\frac{4}{r}(1-r)\log(1-r).
 \end{equation*}
and
\begin{equation*}
	\sum_{n=2}^{\infty}\frac{4(-1)^{n-1}}{n^2+n}=2(-3+4\log 2).
\end{equation*}
On the other hand, we have
\begin{align*}
	\sum_{n=2}^{\infty}\frac{16nr^{2n}}{(n^2+n)^2}&=\sum_{n=2}^{\infty}-\frac{16r^{2n+2}}{r^2(n+1)}+\sum_{n=2}^{\infty}-\frac{16r^{2n}}{(n+1)^2}+\sum_{n=2}^{\infty}\frac{16r^{2n}}{n}\\&=\frac{16}{r^2}\left(\frac{r^4}{2}+r^2+\log(1-r^2)\right)-\frac{16}{r^2}\left(\frac{r^4}{4}-r^2+{\rm Li}_2(r^2)\right)\\&\quad\quad+16\left(-\log(1-r^2)-r^2\right)\\&=-4r^2+\frac{16}{r^2}\left((1-r^2)\log(1-r^2)\right)-\frac{16}{r^2}{\rm Li}_2(r^2)+32.
\end{align*}
Therefore, we have
\begin{align}
	\label{e-33.27}
	&r^2+r+\sum_{n=2}^{\infty}\frac{4r^n}{n^2+n}+\sum_{n=2}^{\infty}\frac{16nr^{2n}}{(n^2+n)^2}\\&=\nonumber	r^2+r+4-2r+\frac{4}{r}(1-r)\log(1-r)-4r^2+\frac{16(1-r^2)\log(1-r^2)}{r^2}\\&\nonumber\quad\quad-\frac{16{\rm Li}_2(r^2)}{r^2}+32.
\end{align}
It is easy to see that
\begin{align*}
&r^2+r+4-2r+\frac{4}{r}(1-r)\log(1-r)-4r^2+\frac{16(1-r^2)\log(1-r^2)}{r^2}\\&\nonumber\quad\quad-\frac{16{\rm Li}_2(r^2)}{r^2}+32\leq 2(-3+4\log 2)
\end{align*}
for $ r\leq r_f(1/2) $, where $ r_f(1/2) $ is root of $ F(r)=0 $ in $ (0,1) $, where $ F : [0,1]\rightarrow\mathbb{R} $ is defined by 
\begin{align*}
	F(r):&=\frac{4}{r}(1-r)\log (1-r)+\frac{16}{r^2}(1-r^2)\log(1-r^2)-\frac{16}{r^2}{\rm Li}_2(r^2)\\&\quad\quad-3r^2-r+29+8\log 2.
\end{align*}
Using the standard argument, we can show that $ F(r) $ has the unqiue root in $ (0,1) $. Let the root be denoted by $ r_f(1/2) $. Further, a simple computation shows that $ r_f(1/2)\approx 0.600881 $. Hence, we have
\begin{align*}
	&\frac{4}{r_f(1/2)}(1-r_f(1/2))\log (1-r_f(1/2))+\frac{16}{r^2_f(1/2)}(1-r^2_f(1/2))\log(1-r^2_f(1/2))\\&\quad\quad-\frac{16}{r^2_f(1/2)}{\rm Li}_2(r^2_f(1/2))-3r^2_f(1/2)-r_f(1/2)+29+8\log 2=0
\end{align*}
which is equivalent to
	\begin{align}\label{e-33.28}
	&r^2_f(1/2)+r_f(1/2)+4-2r_f(1/2)+\frac{4}{r_f(1/2)}(1-r_f(1/2))\log(1-r_f(1/2))-4r^2_f(1/2)\\&\nonumber\quad\quad+\frac{16(1-r^2_f(1/2))\log(1-r^2_f(1/2))}{r^2_f(1/2)}-\frac{16{\rm Li}_2(r^2_f(1/2))}{r^2_f(1/2)}+32\\&\nonumber= 1+2(-3+4\log 2).
\end{align}
To show that $ r_f(1/2)\approx 0.600881  $ is the best possible, we consider the function $ f=f_{1/2} $ defined by \eqref{e-33.7}. A simple computation using \eqref{e-33.8} and \eqref{e-33.26}, for $ f=f_{1/2} $ and $ r=r_f(1/2) $ shows that 
\begin{align*}
	&|z|+\sum_{n=2}^{\infty}\left(|a_n|+|b_n|\right)|z|^n+\frac{S_r}{\pi}\\&=r^2_f(1/2)+r_f(1/2)+4-2r_f(1/2)+\frac{4}{r_f(1/2)}(1-r_f(1/2))\log(1-r_f(1/2))\\&\nonumber\quad\quad-4r^2_f(1/2)+\frac{16(1-r^2_f(1/2))\log(1-r^2_f(1/2))}{r^2_f(1/2)}-\frac{16{\rm Li}_2(r^2_f(1/2))}{r^2_f(1/2)}+32\\&\nonumber= 1+2(-3+4\log 2)\\&=d(f_{1/2}(0), \partial f_{1/2}(\mathbb{D})).
\end{align*}
Therefore, $ r_f(1/2) $ is the best possible. This completes the proof of (i).\vspace{1.5mm}

\noindent (ii) Using Lemma \ref{lem-1.3}, for $ |z|=r $, we obtain
\begin{align}\label{e-33.29}
	&|f(z)|+\sum_{n=2}^{\infty}\left(|a_n|+|b_n|\right)|z|^n+\frac{S_r}{\pi}\\&\nonumber\leq r+\sum_{n=2}^{\infty}\frac{4r^n}{n^2+n}+\left(r^2+\sum_{n=2}^{\infty}\frac{16nr^{2n}}{n^2+n}\right)\\&\nonumber=r+2\left(4-2r+\frac{4}{r}(1-r)\log(1-r)\right)\\&\nonumber\quad\quad+\left(r^2-4r^2+\frac{16(1-r^2)\log(1-r^2)}{r^2}-\frac{16{\rm Li}_2(r^2)}{r^2}+32\right)\\&\nonumber=-3r^2-3r+\frac{8}{r}(1-r)\log(1-r)+\frac{16(1-r^2)\log(1-r^2)}{r^2}-\frac{16{\rm Li}_2(r^2)}{r^2}+40.
\end{align}
A simple computation shows that
\begin{align*}
	&-3r^2-3r+\frac{8}{r}(1-r)\log(1-r)+\frac{16(1-r^2)\log(1-r^2)}{r^2}-\frac{16{\rm Li}_2(r^2)}{r^2}+40\\&\leq 1+2(-3+4\log 2)
\end{align*}
for $ r\leq r^*_f(1/2) $, where $ r^*_f(1/2) $ is a root of $ T(r)=0 $ in $ (0,1) $, where $ T : [0,1]\rightarrow\mathbb{R} $ is defined by
\begin{align*}
	T(r):&=-3r^2-3r+\frac{8}{r}(1-r)\log(1-r)+\frac{16(1-r^2)\log(1-r^2)}{r^2}-\frac{16{\rm Li}_2(r^2)}{r^2}\\&\quad\quad+45-8\log 2.
\end{align*}
By a simple computation, we can show that $ T(r) $ has the unique root in $ (0,1) $. Let $ r^*_f(1/2) $ be the root of $ T(r) $. Then we have $ r^*_f(1/2)\approx 0.302059 $. It is easy to see that
\begin{align*}
	&-3\left(r^*_f(1/2)\right)^2-3r^*_f(1/2)+\frac{8}{r^*_f(1/2)}(1-r^*_f(1/2))\log(1-r^*_f(1/2))\\&\quad\quad+\frac{16(1-\left(r^*_f(1/2)\right)^2)\log(1-\left(r^*_f(1/2)\right)^2)}{\left(r^*_f(1/2)\right)^2}-\frac{16{\rm Li}_2(\left(r^*_f(1/2)\right)^2)}{\left(r^*_f(1/2)\right)^2}\\&\quad\quad+45-8\log 2=0
\end{align*}
which is equivalent to
\begin{align}\label{e-33.30}
	&-3\left(r^*_f(1/2)\right)^2-3r^*_f(1/2)+\frac{8}{r^*_f(1/2)}(1-r^*_f(1/2))\log(1-r^*_f(1/2))\\&\nonumber\quad\quad+\frac{16(1-\left(r^*_f(1/2)\right)^2)\log(1-\left(r^*_f(1/2)\right)^2)}{\left(r^*_f(1/2)\right)^2}-\frac{16{\rm Li}_2(\left(r^*_f(1/2)\right)^2)}{\left(r^*_f(1/2)\right)^2}+40\\&\nonumber=1+2(-3+4\log 2).
\end{align}
To show that $ r^*_f(1/2)\approx 0.302059 $ is the best possible, we consider the function $ f=f_{\alpha} $ defined by \eqref{e-1.4}. Then in view of \eqref{e-3.3}, \eqref{e-33.29} and \eqref{e-33.30}, for $ f=f_{\alpha} $ and $ r= r^*_f(1/2) $, we obtain
\begin{align*}
& |f_{\alpha}(z)|+\sum_{n=2}^{\infty}\left(|a_n|+|b_n|\right)|z|^n+\frac{S_{r^*_f(1/2)}}{\pi}\\&=-3\left(r^*_f(1/2)\right)^2-3r^*_f(1/2)+\frac{8}{r^*_f(1/2)}(1-r^*_f(1/2))\log(1-r^*_f(1/2))\\&\nonumber\quad\quad+\frac{16(1-\left(r^*_f(1/2)\right)^2)\log(1-\left(r^*_f(1/2)\right)^2)}{\left(r^*_f(1/2)\right)^2}-\frac{16{\rm Li}_2(\left(r^*_f(1/2)\right)^2)}{\left(r^*_f(1/2)\right)^2}+40\\&\nonumber=1+2(-3+4\log 2)\\&=d(f_{\alpha}(0), \partial f_{\alpha}(\mathbb{D})).
\end{align*}
This shows that $ r^*_f(1/2) $ is the best possible. This completes the proof of (ii).
\end{proof}	
\begin{proof}[\bf Proof of Theorem \ref{th-2.15}]
	(i) Let $f\in \mathcal{W}^0_{\mathcal{H}}(\alpha) $, then for $ |z|=r $, using Lemma \ref{lem-1.3}, we obtain
	\begin{align}\label{e-3.26}
		|z|+|h(z)|+\sum_{n=2}^{\infty}|a_n||z|^n\leq d(f(0), \partial f(\mathbb{D}))
	\end{align}
if 
\begin{equation*}
	2r+\sum_{n=2}^{\infty}\frac{4r^n}{\alpha n^2+(1-\alpha)n}\leq 1+\sum_{n=2}^{\infty}\frac{2(-1)^{N-1}}{\alpha n^2+(1-\alpha)n}.
\end{equation*}
Let $ J_6 : [0,1]\rightarrow\mathbb{R} $ be defined by
\begin{equation*}
	J_6(r):=2r+\sum_{n=2}^{\infty}\frac{4r^n}{\alpha n^2+(1-\alpha)n}- 1-\sum_{n=2}^{\infty}\frac{2(-1)^{N-1}}{\alpha n^2+(1-\alpha)n}. 
\end{equation*}
It is not difficult to show that $ J_6 $ has the unique root in $ (0,1) $. Let $ r_h(\alpha) $ be the root of $ J_6(r) $ and hence  
\begin{equation}\label{e-3.27}
	2r_h(\alpha)+\sum_{n=2}^{\infty}\frac{4r^n_h(\alpha)}{\alpha n^2+(1-\alpha)n}= 1+\sum_{n=2}^{\infty}\frac{2(-1)^{N-1}}{\alpha n^2+(1-\alpha)n}.
\end{equation}
Now it is enough to show that $ r_h(\alpha) $ is the best possible. To prove this, we consider the function $ f=f_{\alpha}=h_{\alpha}+\overline{g_{\alpha}} $ given by \eqref{e-1.4}.  A simple computation using \eqref{e-3.3}, \eqref{e-3.26} and \eqref{e-3.27}, for $ f=f_{\alpha} $ and $ z=r_h(\alpha) $ shows that
\begin{align*}
	|z|+|h_{\alpha}(z)|+\sum_{n=2}^{\infty}|a_n||z|^n&=	2r_h(\alpha)+\sum_{n=2}^{\infty}\frac{4r^n_h(\alpha)}{\alpha n^2+(1-\alpha)n}\\&= 1+\sum_{n=2}^{\infty}\frac{2(-1)^{n-1}}{\alpha n^2+(1-\alpha)n}\\&=d(f_{\alpha}(0), \partial f_{\alpha}(\mathbb{D})).
\end{align*} 
This shows that $ r_h(\alpha) $ is the best possible. This completes the proof of (i).\vspace{1.5mm}

\noindent (ii) Let $ J_7 : [0,1]\rightarrow\mathbb{R} $ be defined by
\begin{equation*}
	J_7(r):=r+\sum_{n=2}^{\infty}\frac{2r^n}{\alpha n^2+(1-\alpha)n}-1-\sum_{n=2}^{\infty}\frac{2r^{n-1}}{\alpha n^2+(1-\alpha)n}.
\end{equation*}
By a simple calculation, we can show that $ J_7(r) $ has the unique root in $ (0,1) $ which we denote by $ r_g(\alpha) $. Therefore,
\begin{equation}\label{e-3.28}
r_g(\alpha)+\sum_{n=2}^{\infty}\frac{2r^n_g(\alpha)}{\alpha n^2+(1-\alpha)n}=1+\sum_{n=2}^{\infty}\frac{2r^{n-1}}{\alpha n^2+(1-\alpha)n}.
\end{equation}
For $ |z|=r $, using Lemma \ref{lem-1.7}, we obtain
\begin{align}\label{e-3.29}
	|z|+|g(z)|+\sum_{n=2}^{\infty}|b_n||z|^n\leq d(f(0), \partial f(\mathbb{D}))
\end{align}
if 
\begin{equation*}
	r+\sum_{n=2}^{\infty}\frac{2r^n}{\alpha n^2+(1-\alpha)n}\leq 1+\sum_{n=2}^{\infty}\frac{2r^{n-1}}{\alpha n^2+(1-\alpha)n}.
\end{equation*}
for $ r\leq r_g(\alpha) $. In order to show that $ r_g(\alpha) $ is the best possible, we consider the function $ f=f^*_{\alpha} $ defined in \eqref{e-1.8}. In view of \eqref{e-3.1}, it is easy to see that 
\begin{equation}\label{e-3.30}
	d(f^*_{\alpha}(0),\partial f^*_{\alpha}(\mathbb{D}))=1+\sum_{n=2}^{\infty}\frac{2(-1)^{n-1}}{\alpha n^2+(1-\alpha)n}.
\end{equation}
A simple computation using \eqref{e-3.28}, \eqref{e-3.29} and \eqref{e-3.30}, for $ f=f^*_{\alpha} $  and $ z=r_g(\alpha) $ shows that
\begin{align*}
		|z|+|g(z)|+\sum_{n=2}^{\infty}|b_n||z|^n&=r_g(\alpha)+\sum_{n=2}^{\infty}\frac{2r^n_g(\alpha)}{\alpha n^2+(1-\alpha)n}\\&=1+\sum_{n=2}^{\infty}\frac{2r^{n-1}}{\alpha n^2+(1-\alpha)n}\\&=d(f^*_{\alpha}(0),\partial f^*_{\alpha}(\mathbb{D})).
\end{align*}
Therefore, $ r_g(\alpha) $ is the best possible. This completes the proof of (ii).
\end{proof}	
\begin{proof}[\bf Proof of Corollary \ref{cor-2.18}]
	Since $ f\in\mathcal{W}^0_{\mathcal{H}}(\alpha) $, for $ \alpha=1/2 $ and $ |z|=r $, in view of Lemma \ref{lem-1.7}, we obtain
	\begin{align}\label{e-3.36}
		|g(z)|+\sum_{n=2}^{\infty}|b_n||z|^n&\leq 2\sum_{n=2}^{\infty}\frac{2r^n}{n^2+n}\\&\nonumber\leq -2r+4+\frac{4(1-r)\log(1-r)}{r}.
	\end{align}
A simple computation shows that
\begin{equation*}
	-2r+4+\frac{4(1-r)\log(1-r)}{r}\leq 1+(-3+4\log 2).
\end{equation*}
for $ r\leq r^*_g(1/2) $, where $ r^*_g(1/2) $ is a root of $ A(r)=0 $, where $ A : [0,1]\rightarrow\mathbb{R} $ is defined by
\begin{equation*}
	A(r):=-2r+6+\frac{4(1-r)\log(1-r)}{r}-4\log 2.
\end{equation*}
A simple computation shows that $ A(r) $ has the unique root in $ (0,1) $, say $ r^*_g(1/2) $ and hence we have $  r^*_g(1/2)\approx 0.794054 $. Therefore, we have 
\begin{equation*}
	-2r^*_g(1/2)+6+\frac{4(1-r^*_g(1/2))\log(1-r^*_g(1/2))}{r^*_g(1/2)}-4\log 2=0
\end{equation*}
which is equivalent to
\begin{equation}\label{e-3.37}
	-2r^*_g(1/2)+4+\frac{4(1-r^*_g(1/2))\log(1-r^*_g(1/2))}{r^*_g(1/2)}=1+(-3+4\log 2).
\end{equation}
To show that $ r^*_g(1/2) $ is the best possible, we consider the function $ f=f^*_{1/2} $ defined in \eqref{e-1.8}. In view of \eqref{e-3.1}, it is easy to see that 
\begin{equation}\label{e-3.38}
	d(f^*_{1/2}(0),\partial f^*_{1/2}(\mathbb{D}))=1+\sum_{n=2}^{\infty}\frac{2(-1)^{n-1}}{ n^2+n}=1+(-3+4\log 2).
\end{equation}
Then using \eqref{e-3.36}, \eqref{e-3.37} and \eqref{e-3.38}, for $ f=f^*_{1/2} $ and $ r=r^*_g(1/2) $, we obtain
\begin{align*}
		|g(z)|+\sum_{n=2}^{\infty}|b_n||z|^n&=-2r^*_g(1/2)+4+\frac{4(1-r^*_g(1/2))\log(1-r^*_g(1/2))}{r^*_g(1/2)}\\&=1+(-3+4\log 2)\\&=d(f^*_{1/2}(0),\partial f^*_{1/2}(\mathbb{D})).
\end{align*}
Therefore, $ r^*_g(1/2)\approx 0.794054 $ is the best possible.
\end{proof}	
\vspace{2mm}
\noindent\textbf{Acknowledgment:} The first author is supported by the Institute Post Doctoral Fellowship of IIT Bhubaneswar, India, the second author is supported by SERB-CRG.

\end{document}